%% file: QuartCurves.tex
\documentclass[a4paper,10pt]{amsart}
 \usepackage{amsfonts,mathrsfs}
 \usepackage{amsthm}
 \usepackage{amssymb}
 \usepackage{enumerate}
 \usepackage{tikz-cd}
 \usepackage{tikz}
 \usepackage{cleveref}
 \usepackage{url}
 \usepackage{framed}
\theoremstyle{definition}
\newtheorem{Def}{Definition}[section]
\newtheorem{ex}[Def]{Example}
\newtheorem{rem}[Def]{Remark}

\theoremstyle{plain}
\newtheorem{prop}[Def]{Proposition}
\newtheorem{thm}[Def]{Theorem}
\newtheorem{alg}[Def]{Algorithm}
\newtheorem*{thm*}{Theorem}
\newtheorem{lem}[Def]{Lemma}
\newtheorem{cor}[Def]{Corollary}
\newtheorem*{cor*}{Corollary}
\newtheorem{con}[Def]{Conjecture}
\newtheorem*{con*}{Conjecture}
\newtheorem*{frag*}{Question}

\newtheorem*{verm*}{Vermutung}

\include{Operators}
\include{Fonts}

\title[Positive Quadratic Determinantal Representations]{(Positive) Quadratic Determinantal Representations of Quartic Curves and the Robinson Polynomial}
\author{Clemens Brüser}
\address{Technische Universit\"at, Dresden, Germany} 
\email{clemens.brueser@tu-dresden.de}
\author{Mario Kummer}
\email{mario.kummer@tu-dresden.de}
\thanks{Both authors were supported by the DFG grant 502861109.}
\begin{document}


\begin{abstract}
 We prove that every real nonnegative ternary quartic whose complex zero set is smooth can be represented as the determinant of a symmetric matrix with quadratic entries which is everywhere positive semidefinite. We show that the corresponding statement fails for the Robinson polynomial, answering a question by Buckley and {\v{S}}ivic.
\end{abstract}
\maketitle

\section{Introduction}

Let $F\in\R[x_0,x_1,x_2]_{2d}$ be a homogeneous polynomial of degree $2d$. We say that $F$ has a \emph{positive quadratic (determinantal) representation} if there exists a symmetric $d\times d$ matrix $M=M(x_0,x_1,x_2)$ whose entries are homogeneous real polynomials of degree two such that $F=\det(M)$ and $M(a)$ is positive semidefinite for all $a\in\R^3$. Since the determinant of a positive semidefinite matrix is nonnegative, the existence of such a representation implies that $F$ is a \emph{nonnegative polynomial}, meaning that $F(a)\geq0$ for all $a\in\R^3$. Already for $d=3$ it has been open whether conversely every nonnegative $F\in\R[x_0,x_1,x_2]_{6}$ has a positive quadratic representation \cite[\S7]{quarez}. Recently, this problem has also drawn attention in material science \cite[p.~7]{harhov}. Buckley and {\v{S}}ivic conjectured in \cite{buckleysivic} that the answer is negative and they suggest that the so-called \emph{Robinson polynomial}
\begin{equation*}
 R=x_0^6-x_0^4 x_1^2-x_0^2 x_1^4+x_1^6-x_0^4 x_2^2+3 x_0^2 x_1^2 x_2^2-x_1^4 x_2^2-x_0^2 x_2^4-x_1^2 x_2^4+x_2^6.
\end{equation*}
could be such a counterexample. This article establishes two key results. First it establishes that for the Robinson polynomial such a representation does not exist.
\begin{thm}\label{thm:rob}
 The Robinson polynomial $R$ does not have a positive quadratic representation.
\end{thm}

In earlier work by Quarez it has been shown that in degree four, too, there are polynomials that do not admit a positive quadratic representation (see \cite[\S 7]{quarez}). Our second result proves that every such counter-example must necessarily define a singular plane curve.
\begin{thm}\label{thm:quart}
 Every smooth positive ternary quartic has a positive quadratic representation. 
\end{thm}

Throughout this article we will restrict considerations to equivalence classes of determinantal representations.

\begin{Def}
    Let $F = \det(M) = \det(M')$ be two symmetric quadratic representations of $F \in K[x_0, x_1, x_2]$ (with $K \in \{\R, \C\}$). We say that $M$ and $M'$ are \emph{equivalent} quadratic representations, if there exists an invertible matrix of constants $T \in M_d[K]$ satisfying $M' = T^tMT$.
\end{Def}

\section{The Robinson Polynomial}

\subsection{General Observations}

For proving \Cref{thm:rob} we first need to find means to deal with the singularities of the curve defined by the Robinson polynomial. To this end, in this section we lay out preliminary observations and derive necessary conditions for a positive quadratic representation to exist for the Robinson polynomial. We fix a symmetric $d\times d$ matrix $M=M(x_0,x_1,x_2)$ the entries of which  are homogeneous polynomials of degree two such that $M(a)$ is positive semidefinite for all $a\in\R^3$ and such that $F=\det(M)\in\textrm{P}_{2d}$ is not identically zero. We further denote $\pp^2_\R=\Proj(\R[x_0,x_1,x_2])$ and let $X=\cV(F)\subset\pp^2_\R$ be the zero set of $F$, which is a projective scheme of (pure) dimension one over $\R$. First we note that under some conditions one can bound the dimension of the kernel of $M$.

\begin{lem}\label{lem:kersing}
 Let $0\neq a\in\C^3$ such that the kernel of the complex $d\times d$ matrix $M(a)$ has dimension at least two. Then $[a]\in X(\C)$ is a singular point of $X$.
\end{lem}

\begin{proof}
 By the chain rule, the total derivative of $\det(M)$ at $a$ is given by
 \begin{equation*}
     \textrm{d}(\det \circ M)|_a = \textrm{d}\det|_{M(a)} \circ \textrm{d}M|_a.
 \end{equation*}
 Now
 \begin{equation*}
     (\textrm{d} \det|_{M(a)})_{i,j} = \left(\frac{\partial \det}{\partial x_{ij}}\right) M(a) = (-1)^{i+j}\det(M^{(i,j)}(a)),
 \end{equation*}
 where $M^{(i,j)}$ indicates that the $i$-th row and $j$-th column of $M$ have been eliminated, i.e. the $(i,j)$-th cofactor of $M$. By choice of $a$, we know that the determinant of every submatrix of $M(a)$ of size at least $(d-1)\times(d-1)$ equals $0$ and the claim follows.
\end{proof}

\begin{lem}\label{lem:kereal}
 Assume that $d=3$ and that there exists $0\neq a\in\R^3$ such that the kernel of $M(a)$ has dimension at least two. Then $F$ is a sum of squares of homogeneous polynomials of degree three.
\end{lem}

\begin{proof}
 By \cite[Theorem 4.4]{quarez} any positive semidefinite biquadratic form $B(x,y)$ in the variables $x=(x_0,x_1,x_2)^t, y=(y_0,y_1,y_2)^t$ that has an infinite number of zeros in $\pp^2_\R \times \pp^2_\R$ is a sum of squares. We apply this to the biquadratic form given by $B(x,y)=y^tM(x_0,x_1,x_2)y$. By assumption, $M(a)$ is positive semidefinite for arbitrary $a$, hence $B$ is a positive semidefinite form. For our specific choice of $a$, the dimension of the kernel of $M(a)$ being at least two yields an infinite number of zeros. This shows that $B$ is a sum of squares.
 
 We now prove that $\det(M)$ is a sum of squares following ideas in \cite[Proposition 5.6]{quarez}. First note that, since $B$ is a sum of squares, we have $M=A^tA$ for some matrix $A \in \R[x]^{k \times 3}$, where $k \in \N$ and the entries of $A$ are linear homogeneous polynomials in $x$. For $S \in \{1, \dots, k\}$ let $A_S$ (resp. $(A^t)_S$) denote the submatrix $A$ (resp. $A^T$) obtained by deleting all rows (resp. columns) with index not contained in $S$. Using the Cauchy--Binet formula and summing over all $S \subset \{1,\dots, k\}$ with $|S| = 3$, we get
 \begin{align*}
     \det(M) = \sum_{S} \det((A^t)_S) \det(A_S) = \sum_{S} \det((A_S)^t) \det(A_S) = \sum_{S} \det(A_S)^2.
 \end{align*}
 Finally, note that $\det(A_S)$ is homogeneous of degree three, since $|S| = 3$.
\end{proof}

\begin{cor}\label{cor:lb3}
 Assume that $d=3$, that $F$ is not a sum of squares and that all singular points of $X$ are in $X(\R)$. Then the dimension of the kernel of $M(a)$ is at most one for all $0\neq a\in\C^3$.
\end{cor}

\begin{proof}
 This is clear from combining \Cref{lem:kersing} and \Cref{lem:kereal}. 
\end{proof}

In the following, we will associate with $M$ a coherent sheaf on $X$.
The matrix $M$ defines an injective morphism of $\cO_{\pp^2_\R}$-modules:
\begin{equation*}
 M\colon\cO_{\pp^2_\R}(-d-1)^d\to\cO_{\pp^2_\R}(-d+1)^d. 
\end{equation*}
We let $\cF$ be the cokernel of this morphism (of coherent sheaves) so that we obtain a short exact sequence of $\cO_{\pp^2_\R}$-modules:
\begin{equation}\label{eq:res}
 0\to\cO_{\pp^2_\R}(-d-1)^d\to\cO_{\pp^2_\R}(-d+1)^d\to\cF\to0.
\end{equation}
By \cite[Theorem A]{beauville} the sheaf $\cF$ is \emph{arithmetically
Cohen-Macaulay} and its support is $X$. The latter implies that it is the pushforward $\iota_*\cL$ of the coherent sheaf $\cL=\iota^*\cF$ on $X$ where $\iota\colon X\hookrightarrow\pp^2_\R$ is the embedding. The situation is particularly nice when $\cL$ is a line bundle.

\begin{lem}\label{lem:lb2t}
 If $F$ is square-free and the dimension of the kernel of $M(a)$ is at most one for all $0\neq a\in\C^3$, then $\cL$ is a line bundle. If additionally $X$ is integral, then $\cL\otimes\cL\cong\cO_X$.
\end{lem}

\begin{proof}
 We note that $M$ induces a morphism of fibres at closed points $a$ via $M(a)$. For the first part, observe that if $a \in X$ then $\det (M(a)) = 0$. Hence the dimension of the kernel of $M(a)$ is at least one. Together with our assumption, the dimension must be exactly one and therefore the dimension of the cokernel is one, too. Following \cite[Exercise II.5.8]{Hart77}, semi-continuity bounds the dimension of the fibres at all (not necessarily closed) points by $1$. Since the determinant vanishes at all points of $X$, the dimension of the fibres is constant on $X$, and as $F$ is assumed to be square-free, $X$ is reduced and hence $\cF$ is locally free (again by \cite[Exercise II.5.8]{Hart77}). 
 The second part then follows from \cite[Theorem B]{beauville} and the fact that $M$ is symmetric.
\end{proof}

Assume that we are in the situation of \Cref{lem:lb2t}. Given $F$ we can then construct $M$ as follows: We go through all (isomorphy classes of) line bundles $\cL$ on $X$ such that $\cL\otimes\cL\cong\cO_X$ and compute their free resolution. Then we can check whether the free resolution has the correct format and, if yes, whether the occurring matrix is positive semidefinite. For this we need to understand the $2$-torsion points in the group $\Pic(X)$. If $X$ is smooth, then $\Pic(X)\cong\mathbb{T}^g\oplus\Z$ where $\mathbb{T}$ is the circle group and $g=(d-1)(2d-1)$ is the genus of $X$ by \cite[\S2 and \S4]{grossharris}. The zero set of the Robinson polynomial however is not smooth and we will now explain how to understand its Picard group.

\subsection{The Picard Group of the Robinson Polynomial} \label{sec:Rob}
Recall that the Robinson polynomial is defined as
\begin{equation*}
 R=x_0^6-x_0^4 x_1^2-x_0^2 x_1^4+x_1^6-x_0^4 x_2^2+3 x_0^2 x_1^2 x_2^2-x_1^4 x_2^2-x_0^2 x_2^4-x_1^2 x_2^4+x_2^6.
\end{equation*}
Its zero set $X = \cV(R) \subset\pp_\R^2$ has exactly ten singular points $(x_0:x_1:x_2)$, namely those points with homogeneous coordinates in $\{0,1,-1\}$ such that at most one of them is $0$. In particular, all of them are real. We denote this set of points by $Y$. The normalization of $X$ is a real curve of genus $0$ with empty real part and thus isomorphic to $Q=\cV(x^2+y^2+z^2)\subset\pp_\R^2$. The normalization map is given by
\begin{equation*}
 \pi\colon Q\to X,\, (x:y:z)\mapsto((2x^2+z^2)z:x(x^2+2z^2):(x-z)(x+z)y).
\end{equation*}
Using that $R$ is not a sum of squares, by \Cref{cor:lb3} and \Cref{lem:lb2t} any quadratic representation of $R$ must correspond to a $2$-torsion point of $\Pic(X)$. To determine $\Pic(X)$ we consider the injective morphism of sheaves $\cO_X \to \pi_* \cO_Q$ induced by the normalization map $\pi$. We restrict this morphism to multiplicative units, and since $\pi$ is an isomorphism outside the singular locus $Y$ of $X$ we get the following short exact sequence of sheaves on $X$
\begin{equation}\label{eq:shortsequence}
 0 \to \cO_X^\times \to (\pi_*\cO_Q)^\times \to \mathbb{T}_Y \to 0
\end{equation}
where $\mathbb{T}_Y$ is the constant sheaf associated with $\mathbb{T}$ supported on $Y$.

\begin{thm}\label{thm:picrob}
 There is a short exact sequence of abelian groups
 \begin{equation*}
     0 \to \mathbb{T}^{10} \to \Pic(X) \to \Z \to 0.
 \end{equation*}
 In particular the set of 2-torsion points of $\Pic(X)$ is isomorphic to $(\Z/2\Z)^{10}$.
\end{thm}

\begin{proof}
 Applying the long exact sequence in cohomology to \Cref{eq:shortsequence} and using that $\mathbb{T}_Y$ is supported on a scheme of dimension zero yields
 \begin{equation*}
     0 \to \R^\times \to \R^\times \to \mathbb{T}^{10} \to H^1(X, \cO_X^\times) \to H^1(X, \pi_*(\cO_Q)^\times) \to 0.
 \end{equation*}
 Now $\Pic(X)=H^1(X, \cO_X^\times)$, and by \cite[Corollaire 21.8.2]{EGAIV4} we have $\Pic(Q) = H^1(Q, \cO_Q^\times) = H^1(X, \pi_*(\cO_Q)^\times)$. Moreover, by \cite[Tag 0CDY]{stacks-project}) we have an injection $\Pic(Q) \hookrightarrow \Pic(Q_\C) \cong \Z$, hence $\Pic(Q) \cong \Z$. Finally, the map $\R^\times \to \R^\times$ is surjective, so that we get the exact sequence
 \begin{equation*}
     0 \to \mathbb{T}^{10} \to \Pic(X) \to \Z \to 0.
 \end{equation*}
 Every 2-torsion point in $\Pic(X)$ necessarily lies in the kernel of $\Pic(X) \to \Z$, thus must arise as the image of a 2-torsion point in $\mathbb{T}^{10}$. But the 2-torsion subgroup of $\mathbb{T}^{10}$ is (isomorphic to) $(\Z/2\Z)^{10}$.
\end{proof}

It is possible to explicitly compute the image of the morphism $\mathbb{T}^{10} \to \Pic(X)$, i.e., for each $a \in \mathbb{T}^{10}$ we obtain a concrete line bundle $\cL_a$ satisfying (in additive notation) $\cL_a \otimes \cL_b = \cL_{a+b}$ as well as $\cL_0 = \cO_X$. Before describing $\cL_a$ explicitly in \Cref{lem:2torexpl} we fix some notation.

\begin{enumerate}
    \item By $K(Q)$ we denote both the function field of $Q$ as well as the constant sheaf on $Q$ associated with it. We identify the set of rational functions in $K(Q)$ all of whose poles lie in the complement of $U$ with $\cO_Q(U)$ so that $\cO_Q$ becomes a subsheaf of $K(Q)$. 

    \item We write $\pi^{-1}(Y) = \{ q_1, \dots, q_{10} \}$. It is a set of ten closed (but not $\R$-valued) points $q_i \in Q\subset \pp^2_\R$.
\end{enumerate}

Now recall from \Cref{eq:res} (with $d=3$) the short exact sequence
\begin{equation*}
 0\to\cO_{\pp^2_\R}(-4)^3\to\cO_{\pp^2_\R}(-2)^3\to\cL\to0
\end{equation*}
where $\cL$ is a 2-torsion point of $\Pic(X)$.

\begin{lem}\label{lem:globalsec}
 The vector space $H^0(\pp^2_\R, \cL(2))$ has dimension 3 over $\R$.
\end{lem}

\begin{proof}
 Tensoring the exact sequence
 \begin{equation*}
  0\to\cO_{\pp^2_\R}(-4)^3\to\cO_{\pp^2_\R}(-2)^3\to\cL\to0
 \end{equation*}
 with $\cO_{\pp^2_\R}(2)$ yields the short exact sequence
\begin{equation*}
 0\to\cO_{\pp^2_\R}(-2)^3\to\cO_{\pp^2_\R}^3\to\cL(2)\to0.
\end{equation*}
The long exact sequence in cohomology then yields
\begin{equation*}
 0 \to
 H^0(\pp^2_\R, \cO_{\pp^2_\R}(-2)^3) \to H^0(\pp^2_\R, \cO_{\pp^2_\R}^3) \to H^0(\pp^2_\R, \cL(2)) \to
 H^1(\pp^2_\R, \cO_{\pp^2_\R}(-2)^3).
\end{equation*}
Using that both the first and last vector space in this sequence are trivial we get an isomorphism
\begin{equation*}
 \R^3 \cong H^0(\pp^2_\R, \cO_{\pp^2_\R}^3) \cong H^0(\pp^2_\R, \cL(2)).\qedhere
\end{equation*}
\end{proof}

Our goal will be to construct an explicit basis of $H^0(\pp^2_\R, \cL_a(2))$ for $a\in\mathbb{T}^{10}$. For that we will interpret all occurring sheaves as subsheaves of $\pi_*K(Q)$ via the normalization map $\pi$.

\begin{lem}\label{lem:identifyO}
 The sheaf $\cO_X$ is a subsheaf of $\pi_*K(Q)$ and can locally be described by the following data:
 \begin{equation*}
  \cO_X(U) = \left\lbrace f \in\cO_Q(\pi^{-1}(U)) \mid f(q_i) \in \R \text{ for all } i \text{ with } q_i \in \pi^{-1}(U) \right\rbrace.
 \end{equation*}
\end{lem}

\begin{proof}
 The normalization map $\pi$ yields an injective morphism of sheaves
 \begin{equation*}
  0 \to \cO_X \to \pi_*\cO_Q.
 \end{equation*}
 Now $\cO_Q$ is a subsheaf of $K(Q)$ which proves that $\cO_X$ is a subsheaf of $\pi_*K(Q)$.
 
 For the explicit description we apply the results of \Cref{app:realblowup}. Fix a real zero $p$ of $R$ (and thus one of its singularities) and after translation in the real plane assume $p = (0,0,1)$. Let $\hat{F}(x,y)$ denote the resulting polynomial after dehomogenization with respect to $x_2$. One easily verifies by an explicit computation that in our setting the degree two part of $\hat{F}$ is irreducible, so our hypotheses satisfy the conditions of the discussion in \Cref{app:realblowup}. Applying \Cref{prop:Aexpl} to the blowup of $\hat{F}$ at $(x,y)$ yields that the stalk of $\cV(\hat{F})$ in $(x,y)$ is the subset of the stalk of the blow-up of $\cV(\hat{F})$ at $(x,y)$ comprising all those germs whose value in $(x,y)$ lies in $\R$.
\end{proof}

\begin{cor}
 Let $H$ be the divisor of zeros of a linear form $h$ on $X$ such that $\operatorname{supp}(H) \cap Y = \emptyset$. Define
 \begin{equation*}
  h_\pi(x,y,z) = h((2x^2+z^2)z,x(x^2+2z^2),(x-z)(x+z)y)
 \end{equation*}
 and denote its divisor of zeros on $Q$ by $H_\pi$. The sheaf $\cO_X(t)$ is isomorphic to the subsheaf of $\pi_*K(Q)$ that can locally be described by the following data:
 \begin{equation*}
  \cO_X(t)(U) = \left\lbrace f \in \cO_Q(\pi^{-1}(U)\backslash H_\pi) \mid
  \begin{array}{c}
   f(q_i)\in \R \text{ for } i \text{ with } q_i \in \pi^{-1}(U) \\  \text{ and }
   \text{poles in } H_\pi \cap \pi^{-1}(U) \\ \text{ at most of order } t 
  \end{array} \right\rbrace.
 \end{equation*}
\end{cor}

\begin{proof}
 We tensor the sequence
 \begin{equation*}
  0 \to \cO_X \to \pi_*\cO_Q
 \end{equation*}
 by $\cO_X(t)$ and obtain the injection
 \begin{equation*}
  0 \to \cO_X(t) \to \pi_*\cO_Q \otimes \cO_X(t).
 \end{equation*}
 By the projection formula (\cite[Exercise II.5.1.d]{Hart77}) we get an isomorphism
 \begin{equation*}
     \pi_*\cO_Q \otimes \cO_X(t) \cong \pi_*(\cO_Q \otimes \pi^*\cO_X(t)).
 \end{equation*}
 The right hand side is a subsheaf of $\pi_*K(Q)$ by right-exactness of the pushforward functor and because all invertible sheaves on $Q$ are subsheaves of $K(Q)$. The explicit description of $\cO_X(t)$ now follows from the identification of $\cO_X$ in $K(Q)$ given in \Cref{lem:identifyO}.
\end{proof}

\begin{lem}\label{lem:2torexpl}
 Let $a \in \mathbb{T}^{10}$. Then we have the following:
 \begin{enumerate}
     \item The assignment
     \begin{equation*}
         U \mapsto \cL_a(U) := \left\lbrace f \in \cO_Q(\pi^{-1}(U)) \mid f(q_i) \in \R a_i \text{ for all } i \text{ with } q_i \in \pi^{-1}(U) \right\rbrace
     \end{equation*}
     defines a line bundle on $X$.
     \item The map
     \begin{equation*}
        \mathbb{T}^{10} \to \Pic(C), a \mapsto \cL_a
     \end{equation*}
     is an injective homomorphism of groups.
 \end{enumerate}
\end{lem}

\begin{proof}
 We begin by showing that $\cL_a$ as defined is a line bundle. To do so note that by stalk-exactness of \Cref{eq:shortsequence} for each $a \in \mathbb{T}^{10}$ and for each $i$ there exists an open set $U_i \subseteq X$ such that $\pi^{-1}(U_i) \cap Y = \{q_i\}$, and there exists a rational function $f_i \in \cO_Q^\times(\pi^{-1}(U_i))$ satisfying $f(q_i) = a_i$. Then $\cL_a(U_i) = f_i\cO_X(U_i)$ is trivial. Finally, set $U_0 = X \backslash Y$. Then $\cU = (U_i)_{i=0}^{10}$ is an open cover of $X$, which trivializes $\cL_a$.

 It follows from the definition that
 \begin{equation*}
     \mathbb{T}^{10} \to \Pic(X), a \mapsto \cL_a
 \end{equation*}
 is a group homomorphism and we want to prove that it is injective, i.e., no two line bundles of the form $\cL_a, \cL_b$ are isomorphic for $a \neq b$. To see this assume that $\cL_a$ has local basis $(f_i)_i$ (as constructed above) and $\cL_b$ has local basis $(g_i)_i$. If $\cL_a \cong \cL_b$, then the functions $f_i$ and $g_i$ have to agree on their common domain (containing $q_i$) - up to multiplication with elements in $K(X)^\times$. By \Cref{lem:identifyO} this contradicts the choices of $f_i$ and $g_i$.
\end{proof}

\begin{cor}
 Every 2-torsion point of $\Pic(X)$ is of the form $\cL_a$ for some $a \in \mathbb{T}^ {10}$.
\end{cor}

\begin{proof}
 We already know by \Cref{thm:picrob} that the number of 2-torsion points in $\Pic(X)$ is $2^{10}$. Furthermore, $\cL_a$ is a 2-torsion point for each $a \in T^{10} \cap \{1,i\}^{10}$. By \Cref{lem:2torexpl} the assignment $a \mapsto \cL_a$ is injective, which proves the claim.
\end{proof}

\begin{rem}
 One can in fact show that the homomorphism $\mathbb{T}^{10} \to \Pic(X)$ constructed in \Cref{lem:2torexpl} coincides with that induced by taking cohomology in \Cref{thm:picrob}.
\end{rem}

\subsection*{Representing Types}

The approach to finding determinantal representations as laid out above follows the ideas in \cite[Theorem B]{beauville}. However, not every 2-torsion point gives rise to a symmetric quadratic representation. For necessary and sufficient conditions on the 2-torsion point $\cL$ for the resulting determinantal representation to be quadratic, see \cite[Theorem 6.2.1]{dolgachev} or \cite[Proposition 3.5]{beauville}. These conditions not being satisfied do not present an obstacle to finding a determinantal representation of $R$, though. It just might (and does) happen that the distribution of degrees across the representing matrix is different.

\begin{ex}
 The following two matrices are determinantal representations of the Robinson polynomial:
 \begin{align*}
      & \begin{pmatrix}
         4x_0^2-x_1^2+4x_2^2 & -6x_0^2-2x_0x_1+4x_1^2-6x_2^2 & 4x_0x_1+x_1^2 \\
         & 4x_0^2+8x_0x_1-8x_1^2+12x_2^2 & -2x_0^2-4x_0x_1-2x_1^2 \\
         & & x_0^2+2x_0x_1+4x_1^2-3x_2^2
     \end{pmatrix}, \\
      & \begin{pmatrix}
         x_0^4+6x_0^2x_1^2+5x_1^4-2x_0^2x_2^2-7x_1^2x_2^2+x_2^4 & -x_0^3-x_0x_1^2+x_0x_2^2 & x_0^2+x_1^2-x_2^2 \\
        & -x_1^2-x_2^2 & -x_0 \\
        & & 0
     \end{pmatrix}.
 \end{align*}
 The first representation corresponds to the tuple $(i,i,1,1,i,i,1,1,1,1)$. The second one, corresponding to the tuple $ (1, 1, 1, 1, i, i, 1, 1, 1, 1)$, does not have only quadratic entries.
\end{ex}

Concretely, in the case of the Robinson polynomial $R$ we obtain the following survey of \emph{representing types} of $R$. Here, each entry in a matrix represents the degree of the polynomial in that position, with $-\infty$ indicating the constant 0 polynomial. For convenience, we write $d^{k\times l}$ to denote a $k \times l$-matrix, the entries of which are all of degree $d$.
\begin{align*}
    M_3 &:= (6) \\
    M_2 &:= \begin{pmatrix}
    4 & 4 & 1 \\
    4 & 4 & 1 \\
    1 & 1 & -\infty
    \end{pmatrix} \\
    M_1^k &:= \begin{pmatrix}
    4 & 3 & 2^{1\times k} \\
    3 & 2 & 1^{1\times k} \\
    2^{k\times1} & 1^{k\times1} & -\infty^{k\times k}
    \end{pmatrix}, &k = 0, 1 \\
    M_0^k &:= \begin{pmatrix}
    2^{3\times3} & 1^{3 \times k} \\
    1^{k \times 3} & -\infty^{k \times k}
    \end{pmatrix}, &k = 0, 1, 2
\end{align*}

This gives rise to the following table of dimensions. In particular when going through the 2-torsion points of $\Pic(X)$ with the sole interest of identifying (positive) quadratic representations,  it may first be checked if $H^0(X, \cL(1)) = 0$.

\begin{table}[h]
    \centering
    \begin{tabular}{c|c|c|c|c}
        & $\dim(H^0(\cL_a(0)))$ & $\dim(H^0(\cL_a(1)))$ & $\dim(H^0(\cL_a(2)))$ & $\dim(H^0(\cL_a(3)))$ \\ \hline
        $M_3$ & 1 & 3 & 6 & 10\\
        $M_2$ & 0 & 2 & 5 & 9\\
        $M_1^k$ & 0 & 1 & 4 & 9 \\
        $M_0^k$ & 0 & 0 & 3 & 9
    \end{tabular}
    \caption{Dimension of twists of 2-torsion points.}
    \label{tab:twistcohom}
\end{table}

\subsection*{Explicit results}

In \Cref{app:RobComp} we give a walkthrough of the code producing the free resolutions of all 1024 different 2-torsion points. The result yields the following classification of representing types.

\begin{table}[h!]
   \centering
   \begin{tabular}{c|c|c}
       type &  nr. of occurrences & example for $a\in\mathbb{T}^{10}$ \\ \hline
       $M_3$ & 1 & $(1,1,1,1,1,1,1,1,1,1)$ \\
       $M_2$ & 10 & $(i,1,1,1,1,1,1,1,1,1)$ \\
       $M_1^0$ & 42 & $(i,1,1,i,i,1,1,1,1,1)$ \\
       $M_1^1$ & 45 & $(i,i,1,1,1,1,1,1,1,1)$ \\
       $M_0^0$ & 438 & $(i,i,i,i,1,1,1,1,1,1)$ \\
       $M_0^1$ & 378 & $(i,1,i,1,i,1,1,1,1,1)$ \\
       $M_0^2$ & 110 & $(i,i,i,1,1,1,1,1,1,1)$
    \end{tabular}
    \caption{Representing types of $R$.}
    \label{tab:classesRob}
\end{table}

The matrices of interest for quadratic represenations are those of type $M_0^0$. To decide if these give a positive quadratic representation, it suffices to evaluate them in an arbitrary real point, which the Robinson polynomial does not vanish in, and check if the resulting real-valued matrix is positive definite. Due to the continuity of eigenvalues, the locus in which the representation is positive definite is open. Similarly, the set of all points of evaluation in which the representation is invertible and has at least one negative eigenvalue is open. These two sets cover all of $\pp^2_\R(\R) \backslash Y$ which is connected since $Y$ is finite. This implies that one of the sets must be empty.

Finally, the concrete calculations yield that no quadratic representation of type $M_0^0$ of the Robinson polynomial $R$ is positive. This proves \Cref{thm:rob}.

\begin{rem}
 The methods that we have employed to derive our result about the Robinson polynomial are applicable more generally. For example, whenever we are given a nonnegative sextic curve with ten real zeros of multiplicity one, then its normalization is isomorphic to $\cV(x^2 + y^2 + z^2)$, and we may proceed as for the Robinson polynomial. 
\end{rem}

\section{Smooth Positive Ternary Quartics}

\subsection{The Intrinsic Geometry of a Curve}

We first recall some of the intrinsic geometry underlying the set of all quadratic representations of a given smooth projective \emph{complex} curve $C = \cV(F)$ of degree $\deg(F) = 2d$ and genus $g = (d-1)(2d-1)$. We denote by $J := J(C)$ its Jacobian variety (i.e. the group scheme of degree zero divisor classes on $C$). The foundation, from which our theory is developed lies in the following classical statement, which follows from \cite[Theorem B]{beauville}. It is also established in \cite[\S6.2.2]{dolgachev}.

\begin{thm} \label{thm:2tor=QDR}
    There is a bijection between invertible sheaves $\cL$ of degree zero on $C$ satisfying
    \begin{equation*}
        h^0(C, \cL(d-2)) = 0 = h^1(C, \cL(d-1))
    \end{equation*}
    and equivalence classes of quadratic representations of $F$ over $\C$.
\end{thm}

\begin{rem}
    If we restrict considerations to symmetric determinantal representations, then there is a bijection to the 2-torsion points $\cL$ of the group $J(\C)$ of the Jacobian $J$ with the sole condition
    \begin{equation*}
        h^0(C, \cL(d-2)) = 0.
    \end{equation*}
    This uses $\omega_C \cong \cO_C(2d-3)$ and Serre duality.
\end{rem}

The number of nontrivial 2-torsion points of $J(\C)$ is $2^{2g}-1$, and by \cite[Proposition 4.6]{beauville} for a general curve every nontrivial 2-torsion point $\cL$ of $J(\C)$ will satisfy these conditions. We will next describe how 2-torsion points can be understood through theta characteristics on $C$.

\begin{Def}
    Let $C$ be a smooth plane projective curve of degree $2d$ over $\C$. A line bundle $\theta$ on $C$ is called a \emph{theta characteristic}, if $\theta^{\otimes2} \cong \omega_C \cong \cO_C(2d-3)$. We call a theta characteristic $\theta$ even (resp. odd), if $h^0(C, \theta)$ is even (resp. odd). We call a theta characteristic $\theta$ \emph{vanishing}, if $h^0(C, \theta) > 1$.
\end{Def}

\begin{Def}
    Let $C$ and $C'$ be projective plane curves over $\C$. We say that $C'$ is a \emph{contact curve} to $C$, if the intersection divisor
    \begin{equation*}
        C.C' := \sum_{P \in C} I(P; C,C') \in \Div(C)
    \end{equation*}
    is divisible by 2, i.e., there exists $D \in \Div(C)$ satisfying $C.C' = 2D$. We call $D$ the \emph{contact divisor} of $C$ and $C'$. If $C = \cV(F), C' = \cV(G)$ are reduced curves, then we will also write $F.G = C.C'$.
\end{Def}

Geometrically, odd theta characteristics characterize contact curves of degree $2d-3$ to $C$ (unique ones, if the theta characteristic is non-vanishing). More precisely, let $D$ be an effective divisor on $C$ representing an odd theta characteristic $\theta_-$. Then $D$ is a sum of points on $C$ such that $2D \sim (2d-3)H$, the divisor of zeros on $C$ of a curve of degree $2d-3$.

Non-vanishing even theta characteristics $\theta_+$, on the other hand, give rise to contact curves of degree $2d-1$ (via effective representatives of the divisor class $[T]$ satisfying $\theta_+(1) \cong \cO(T)$). These theta characteristics are in bijection with linear determinantal representations of $C$ (see, for example, \cite[Theorem 2, Theorem 4, Theorem 5]{Vinnikov1989LDR}).

One important result linking odd theta characteristics and 2-torsion points of $J(\C)$ is the following.

\begin{lem} \label{lem:2tor=DiffOddTheta}
    Let $\theta_1, \theta_2$ be two odd theta characteristics. Then $\cL = \theta_1 \otimes \theta_2^\vee$ is a 2-torsion point of $J(C)$. Conversely, every 2-torsion point $\cL$ of $J(C)$ can be expressed as $\cL = \theta_1 \otimes \theta_2^\vee$ with $\theta_1, \theta_2$ being odd theta characteristics.
\end{lem}

\begin{proof}
    The first statement is clear. For the converse, see the proof of \cite[Proposition 5.4.5]{dolgachev}.
\end{proof}

Hence, one way to understand the quadratic representations of $F$ is via contact curves of degree $2d-3$ to $C = \cV(F)$.

A different approach links quadratic representations of $F$ to contact curves of degree $2d-2$. Following \cite[\S 6.2]{dolgachev}, given a 2-torsion point satisfying the conditions of \Cref{thm:2tor=QDR}, a system of contact curves of degree $2d-2$ is given by effective representatives of the unique divisor class $[A]$ satisfying $\cL(d-1) \cong \cO(A)$. Precisely, one has the following lemma (\cite[p.~236--237]{dolgachev}, but see also \cite[Proposition 5.7]{plaumann} for the special case of $d =2$ and \cite[Proposition 3.2]{hvelemt} for a similar result).

\begin{lem}
    Let $Q$ be a quadratic representation of $F$. Then the set
    \begin{equation*}
        \cC_Q := \left\lbrace \lambda^t Q^{\adj} \lambda \mid \lambda \in \pp^{d-1}_\C \right\rbrace
    \end{equation*}
    is a system of contact curves of degree $2d-2$ to $C$. If $a, b \in \cC_Q$, then their divisors of zeros on $C = \cV(F)$ are linearly equivalent.
\end{lem}

Conversely, if we know at least one such contact curve $a \in \cC_Q$, then there is an algorithm that recovers the quadratic representation $Q$ of $F$. This is an adaptation of Dixon's algorithm (see \cite{Dixon1902LDR} for the first description and \cite{Vinnikov1989LDR} for a more modern account), which was originally designed for the construction of linear determinantal representations. An approach very similar to ours has been given in \cite[\S 4]{hvelemt} for linear determinantal representations of hyperbolic plane curves. One of its key components is the following theorem.

\begin{thm}[Max Noether's Fundamental Theorem] \label{thm:AF+BG}
    Let $F,G,H$ be homogeneous polynomials such that $\cV(F)$ and $\cV(G)$ have no common component. Furthermore assume that $F.H \geq F.G$ as divisors on $\cV(F)$. Then there exist homogeneous polynomials $p,q$ satisfying
    \begin{equation*}
        H = pF + qG.
    \end{equation*}
    If $F,G,H$ are all real curves, then $p,q$ may be chosen real as well.
\end{thm}

\begin{proof}
    See \cite[\S 5.5]{FultonCurves} and \cite[Theorem 4.1]{hvelemt}.
\end{proof}

\begin{alg}[Dixon's Algorithm for Quadratic Representations of Ternary Quartics] \label{alg:Dixon}
Let $C = \cV(F)$ be a smooth projective plane curve of degree $4$ over $\C$ and let $C' = V(a_{00})$ be a reduced contact curve to $C$ of degree two with contact divisor $A$. Then the following algorithm returns a quadratic representation $Q$ of $F$, such that $a_{00} \in \cC_Q$.
\begin{enumerate}
    \item Choose a conic $a_{10} \in \cL(2H - A)$ such that the pair $a_{00}, a_{10}$ is linearly independent.
    \item From Max Noether's fundamental theorem choose a conic $a_{11}$ and $p \in \C$ satisfying
    \begin{equation*}
        a_{00}a_{11} - a_{10}^2 = pF.
    \end{equation*}
    \item Define $Q$ as the symmetric matrix
    \begin{equation*}
        Q = \begin{pmatrix}
            a_{00} & a_{10} \\
            a_{10} & a_{11}
        \end{pmatrix}.
    \end{equation*}
\end{enumerate}
\end{alg}

\begin{proof}
First observe that by Riemann--Roch we have
\begin{equation*}
    \ell(2H - A) \geq \deg(2H - A) + 1 - g_C = 4 + 1 - 3 = 2.
\end{equation*}
In particular we may choose $a_{10}$ as desired. Next we know that
\begin{equation*}
    a_{10}^2.F \geq a_{00}.F
\end{equation*}
by construction, so the hypotheses of Max Noether's fundamental theorem (\Cref{thm:AF+BG}) are satisfied. Furthermore, note that for degree reasons we have $p \in \C$, so unless $p = 0$ we have found a quadratic representation.

We argue that $p = 0$ if and only if $a = a_{00}$ is non-reduced. Since $d = 2$ this is the case if and only if $a = l^2$ is the square of a linear form. Indeed, if $p = 0$, then $a_{00}a_{11} = a_{10}^2$. By unique factorization, and since $a_{00},a_{10}$ are linearly independent, $a_{10}^2 = l_1^2l_2^2$ for certain linear factors $l_1,l_2$. Again by linear independence without loss of generality we get $l_1^2 = a_{00}$. Conversely, assume $a = l_1^2$. Then $a_{10}$ vanishes at 4 colinear points, hence is a reducible conic $a_{10} = l_1l_2$ itself. This implies $a_{11} = l_2^2$, so $p = 0$.
\end{proof}

\begin{rem}
The proof of \Cref{alg:Dixon} is elementary in the case of quartic curves. A similar algorithm can be used to construct quadratic representations of arbitrary smooth curves of even degree. The proof relies on the fact that if $a_{00}$ is a reduced contact curve to $C$ with contact divisor $A$, then $A - (d-1)H$ defines a 2-torsion point, satisfying the cohomological hypotheses of \Cref{thm:2tor=QDR}. Hence, a determinantal representation is implicitly already contained in the divisor $A$.
\end{rem}

\begin{rem} \label{rem:realDixon}
While \Cref{alg:Dixon} is phrased for a complex curve $C$, it can be adapted to the real setting: For a real quadratic representation of a real smooth projective curve $C = \cV(F)$ with (real) Jacobian $J$ to exist it is necessary that the corresponding 2-torsion point of $J(\C)$ has a representing divisor that is invariant under complex conjugation. Since Max Noether's fundamental theorem holds over $\R$, Dixon's algorithm may still be carried out. 
\end{rem}

While a 2-torsion point as in \Cref{rem:realDixon} could arguably be called a ``real'' 2-torsion point, there is a subtlety not to be overlooked. Namely, there are two notions of reality for 2-torsion points (and more generally divisor classes).

\begin{Def} \label{def:realDivClass}
    Let $C$ be a smooth projective curve over $\R$, and let $D$ be a divisor on $C_\C$.
    \begin{enumerate}
        \item If $[D] = [\overline{D}]$, then $[D]$ is called a \emph{real} divisor class.
        \item If there exists $D_0 \sim D$ such that $\overline{D_0} = D_0$, then $[D]$ is called a \emph{strongly real} divisor class.
    \end{enumerate}
\end{Def}

Note that if $D \sim D'$ and $[D] = [\overline{D}]$, then $[D - D'] = [0] = [\overline{D - D'}]$, so $[D'] = [\overline{D']}$. Hence reality is independent of the chosen representative of $[D]$.

Real divisor classes correspond to real closed points on the Jacobian of $C$ (considered as a scheme over $\R$), and alternatively they can be identified with line bundles in $\Pic(C_\C)$ that are isomorphic to their complex conjugate. Strongly real divisor classes correspond to line bundles in $\Pic(C)$.

If we only consider 2-torsion points, then those used in \Cref{rem:realDixon} are the strongly real ones. Obviously, strongly real implies real, and for 2-torsion points the two notions coincide if and only if the genus of $C$ is even (see \cite[Proposition 3.3]{grossharris}).

\subsection{Quartic Curves: Proof of Existence}

We now restrict our considerations to $C = V(F)$ being a smooth ternary quartic curve defined over $\R$ with $F > 0$. We let $J$ be its Jacobian variety (considered as an $\R$-scheme). There are exactly 63 theta characteristics on $C$, 35 are even and 28 are odd. They are all non-vanishing, and in particular there are exactly 28 bitangent lines (exactly four of which are real). We first prove that the conditions in \Cref{thm:2tor=QDR} are automatically satisfied.

\begin{lem}[{\cite[Theorem 6.2.3]{dolgachev}}]
    Every non-trivial 2-torsion point of $J(\C)$ defines a (complex) symmetric quadratic representation of $F$.
\end{lem}

\begin{proof}
    Let $D$ be a 2-torsion point of $J(\C)$. Then it may be realized as the difference of two odd theta characteristics, which we identify as a difference of divisors $B_1 - B_2$ by \Cref{lem:2tor=DiffOddTheta}. Following \Cref{thm:2tor=QDR} it now suffices to prove
    \begin{equation*}
        \ell(D) = 0.
    \end{equation*}
    Indeed, if $\ell(D) > 0$, then $D$ is a special divisor, so by Clifford's theorem (\cite[Theorem IV.5.2]{Hart77}) we get $2(\ell(D)-1) \leq \deg(D) = 0$. Thus, $\ell(D) = 1$ and $D$ must be the zero divisor, a contradiction to $D$ being non-trivial.
\end{proof}

A fact that we will use repeatedly is the following.

\begin{prop} \label{prop:realsymmLDR}
    There exist real symmetric linear determinantal representations of $F$, and they are in bijection with non-vanishing strongly real even theta characteristics on $C$.
\end{prop}

\begin{proof}
    See \cite[p.~456]{vin93}.
\end{proof}

For a very accessible classification of all 2-torsion points we recall the definition of a Steiner complex.

\begin{Def}
    Let $[D]$ be a 2-torsion point of $J(\C)$. We call the set $\Sigma(D)$ of 2-sets of odd theta characteristics $\{B_1, B_2\}$ satisfying $[B_1 - B_2] = [D]$ the \emph{Steiner complex} associated with $[D]$.
\end{Def}

By abuse of terminology we will sometimes refer to elements of $\Sigma(D)$ as pairs of odd theta characteristics (or pairs of bitangent lines). Since $B_1 - B_2 \sim B_2 - B_1$ and since every bitangent line is uniquely identified with its contact divisor on $C$, there is no significant harm in doing so. The enumerative and combinatorial structure of Steiner complexes on a smooth plane curve $C$ over $\C$ is well understood, see for example \cite[\S 5.4.2]{dolgachev}. For ternary quartics specifically there is a variety of different characterizations of Steiner complexes. To cite them we first define the bitangent matrix of the curve $C$ (as introduced by \cite[\S 3]{plaumann}). It is based on a chosen linear determinantal representation $M$ of $F$.

\begin{prop} \label{prop:bitmatrix}
    Let $M = x_0M_0 + x_1M_1 + x_2M_2$ be a (complex) symmetric determinantal representation of $F$. Then the three quadric surfaces $\cV(u^t M_i u) \subseteq \pp^3_\C$ intersect in 8 distinct points $O_i$. Denote by $O$ the matrix with rows $O_1, \dots, O_8$. Then
    \begin{equation*}
        \cB = O M O^t
    \end{equation*}
    is a symmetric matrix whose diagonal entries are zero and its other entries are the linear forms defining the 28 bitangent lines of $C$.
\end{prop}

\begin{proof}
    See \cite[\S 6.3.2]{dolgachev} and \cite[Proposition 3.3]{plaumann}.
\end{proof}

The collection of points $O_i$ from \Cref{prop:bitmatrix} is called the \emph{Cayley octad} associated with $M$, the matrix $\cB$ is called the \emph{bitangent matrix} associated with $M$. It has the form
\begin{equation*}
    \cB = \begin{pmatrix}
        0 & b_{12} & b_{13} & b_{14} & b_{15} & b_{16} & b_{17} & b_{18} \\
        b_{12} & 0 & b_{23} & b_{24} & b_{25} & b_{26} & b_{27} & b_{28} \\
        b_{13} & b_{23} & 0 & b_{34} & b_{35} & b_{36} & b_{37} & b_{38} \\
        b_{14} & b_{24} & b_{34} & 0 & b_{45} & b_{46} & b_{47} & b_{48} \\
        b_{15} & b_{25} & b_{35} & b_{45} & 0 & b_{56} & b_{57} & b_{58} \\
        b_{16} & b_{26} & b_{36} & b_{46} & b_{56} & 0 &  b_{67} & b_{68} \\
        b_{17} & b_{27} & b_{37} & b_{47} & b_{57} & b_{67} & 0 & b_{78} \\
        b_{18} &  b_{28} & b_{38} & b_{48} & b_{58} & b_{68} & b_{78} & 0
    \end{pmatrix}
\end{equation*}
and it satisfies $b_{ij} = O_iMO_j^T$ for all $i,j$.

It is now possible to state a concise characterization of all Steiner complexes of a plane quartic. This is a classical result, the formulation of which we have taken from \cite[Theorem 5.10]{plaumann}.

\begin{thm} \label{thm:steiner}
    Let $S = \{ \{l_1, l_1'\}, \dots, \{l_6, l_6'\} \}$ be a set of six pairs of bitangent lines. Then the following statements are equivalent:
    \begin{enumerate}
        \item Let $B_i$ and $B_i'$ denote the contact divisor on $C$ of $\cV(l_i)$ and $\cV(l_i')$ respectively. Then $B_i - B_i' \sim B_j - B_j'$ for all $1 \leq i,j \leq 6$, and the Steiner complex associated with their class is exactly $S$.
        
        \item There exists a quadratic representation $Q$ of $F$ such that $l_il_i' \in \cC_Q$ for each $1 \leq i \leq 6$.

        \item For every choice $1 \leq i,j \leq 6, i \neq j$ the eight points where $\cV(l_il_i'l_jl_j')$ meets $C$ lie on a conic.

        \item Write $[8] := \{1, \dots, 8\}$. Then $S$ is of one of the following forms:
        \begin{align*}
            S^{(1)}_I = \left\lbrace \{b_{ik}, b_{jk} \} \mid \{i,j\} = I, k \in I^c \right\rbrace, \qquad & I \in \binom{[8]}{2} \\
            S^{(2)}_I = \left\lbrace \{b_{ij}, b_{kl} \} \mid \{i,j,k,l\} \in \{I, I^c\} \right\rbrace, \qquad & I \in \binom{[8]}{4}
        \end{align*}
    \end{enumerate}
\end{thm}

\begin{proof}
    See \cite[Theorem 5.10]{plaumann} for references.
\end{proof}

We say that a Steiner complex is of type 1 (resp. 2), if it is of the form $S^{(1)}_I$ (resp. $S^{(2)}_I$). Each Steiner complex $S$ may be represented visually as follows. For each pair of bitangent lines contained in $S$ we draw a line connecting the corresponding entries in the bitangent matrix (for reasons of symmetry we only draw lines in either the upper or the lower triangle of the matrix). In \Cref{fig:Steiner} the orange and purple lines indicate Steiner complexes of type 2 with $I = \{1,2,3,4\}$ and $I = \{1,3,6,7\}$ respectively. The blue lines indicate a Steiner complex of type 1 with $I = \{3,4\}$.

\begin{center}
\begin{figure}[h!]
\begin{tikzpicture}
\draw[step=0.5cm,gray,very thin] (0,0) grid (4,4);
\fill[black] (0.5,3.5) rectangle (0,4);
\fill[black] (0.5,3.5) rectangle (1,3);
\fill[black] (1.5,2.5) rectangle (1,3);
\fill[black] (1.5,2.5) rectangle (2,2);
\fill[black] (2.5,1.5) rectangle (2,2);
\fill[black] (2.5,1.5) rectangle (3,1);
\fill[black] (3.5,0.5) rectangle (3,1);
\fill[black] (3.5,0.5) rectangle (4,0);

\draw[orange, thick] (.75,3.75) -- (1.75,2.75);
\draw[orange, thick] (1.75,3.75) -- (1.25,3.25);
\draw[orange, thick] (1.25,3.75) -- (1.75,3.25);
\draw[orange, thick] (2.75,1.75) -- (3.75,0.75);
\draw[orange, thick] (3.75,1.75) -- (3.25,1.25);
\draw[orange, thick] (3.25,1.75) -- (3.75,1.25);

\draw[blue, thick] (1.25,3.75) -- (1.75,3.75);
\draw[blue, thick] (1.25,3.25) -- (1.75,3.25);
\draw[blue, thick] (2.25,2.75) -- (2.25,2.25);
\draw[blue, thick] (2.75,2.75) -- (2.75,2.25);
\draw[blue, thick] (3.25,2.75) -- (3.25,2.25);
\draw[blue, thick] (3.75,2.75) -- (3.75,2.25);

\draw[purple, thick] (0.25,0.75) -- (1.25,1.25);
\draw[purple, thick] (0.75,0.25) -- (1.75,1.75);
\draw[purple, thick] (0.25,1.25) -- (1.25,0.75);
\draw[purple, thick] (0.75,1.75) -- (1.75,0.25);
\draw[purple, thick] (0.25,2.75) -- (2.75,0.75);
\draw[purple, thick] (0.75,2.25) -- (2.25,0.25);
\end{tikzpicture}
\caption{Three real Steiner complexes.} \label{fig:Steiner}
\end{figure}
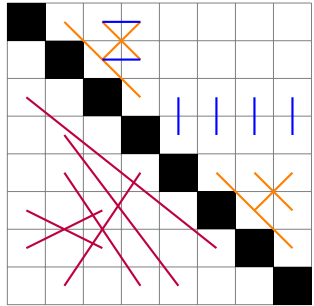
\end{center}

From now on we will assume that $F$ is positive and that the linear determinantal representation $M$ is real. Such a representation $M$ exists by \Cref{prop:realsymmLDR}. In this case the Cayley octad comprises four pairs of complex conjugate points since otherwise $C$ had more than four real bitangent lines. Thus, from now on, in the notation of \Cref{prop:bitmatrix}, we assume without loss of generality $O_{2i} = \overline{O_{2i-1}}$ for $i = 1, \dots, 4$. This implies certain relations via complex conjugation in the bitangent matrix $\cB$, e.g. $b_{13} = \overline{b_{24}}$ or $b_{36} = \overline{b_{45}}$. In particular, we identify the four real bitangent lines of $C$ as $b_{12}, b_{34}, b_{56}$, and $b_{78}$.

\begin{prop}
    Let $\sigma$ denote the permutation $(12)(34)(56)(78) \in S_8$, let $D$ be a divisor whose class $[D]$ is a 2-torsion point, and let $\Sigma(D)$ be the associated Steiner complex. Then $[D]$ is real (in the sense of \Cref{def:realDivClass}) if and only if one of the following holds:
    \begin{enumerate}
        \item $\Sigma(D) = S^{(1)}_I$ and $\sigma(I) = I$.
        \item $\Sigma(D) = S^{(2)}_I$ and $\sigma(I) \in \{I, I^c\}$.
    \end{enumerate}
    It is strongly real if and only if in either case $\sigma(I) = I$.
\end{prop}

\begin{proof}
    First let $\Sigma(D)$ be of type 1, so $D = B_{ik} - B_{jk}$ for $I = \{i, j\}$. Then $D$ is real if and only if $D \sim \overline{D}$, or equivalently $B_{ik} - B_{jk} \sim B_{\sigma(i)\sigma(k)} - B_{\sigma(j)\sigma(k)}$. This is the case if and only if $\sigma(I) = I$.
    
    Next let $\Sigma(D)$ be of type 2. Then the class of $D \sim B_{ij} - B_{kl}$ is real if and only if $B_{ij} - B_{kl} \sim B_{\sigma(i)\sigma(j)} - B_{\sigma(k)\sigma(l)}$. This is the case if and only if $\sigma(I) \in \{I, I^c \}$.

    For the second part of the statement we use that it follows from \cite[p.~727--728]{plaumann} that the quadratic representation $Q$ associated with $D$ is real (or rather is equivalent to one such representation), if and only if $\cC_Q$ contains real contact conics, which is the case if and only if $D$ is such that $\sigma(I) = I$.
\end{proof}

In this setting a 2-torsion point $[D]$ is real if and only if the visual representation of the associated Steiner complex $\Sigma(D)$ is invariant under the simultaneous action of $\sigma$ on the rows and columns of $\cB$.
The proposition furthermore shows that there are exactly seven Steiner complexes that give rise to real quadratic representations:
\begin{enumerate}
    \item The four Steiner complexes of type 1 corresponding to the index sets $I =\{2i-1, 2i\}$ for $1 \leq i \leq 4$.
    \item The three Steiner complexes of type 2 which are associated with $I = \{1,2,3,4\}$, $I = \{1,2,5,6\}$ or $I = \{1,2,7,8\}$.
\end{enumerate}

Finally, we have the necessary tools to prove \Cref{thm:quart}. We do so in a sequence of lemmata. First recall our standing assumptions about $B, T, A$: Both $B$ and $T$ are defined over $\R$, and we define $0 \leq A \sim B + T - H$. Then we define $Q$ as the unique quadratic representation of $F$ corresponding to the divisor class $[D] = [A - H]$. Since $[A-H]$ is defined over $\R$, the system of contact conics $\cC_Q$ contains conics defined over $\R$. By $\cC_Q(\R)$ we denote the set of all real contact conics to $C$.

\begin{lem} \label{lem:real_azyg}
    If $B_1 - B_2$ is the difference of a pair of complex conjugate bitangent lines (identified with their contact divisors), then there exist two real bitangent lines $B_3, B_4$ such that $B_1 - B_2 \sim B_3 - B_4$. There are four 2-torsion points that arise in this way (among them the trivial point that is the class of the zero divisor). In particular, the product $l_1l_2$ lies in $\cC_Q$ if and only if $l_3l_4$ does.
\end{lem}

\begin{proof}
    The trivial 2-torsion point arises from $B_1$ being a real bitangent line. Henceforth assume that this is not the case.
    
    We earlier noted that the real bitangent lines are $b_{12}, b_{34}, b_{56}$, and $b_{78}$. Combining any two of them we retrieve exactly the three index sets (or their complements) corresponding to the three non-trivial real Steiner complexes of type 2.

    If now $b_{ij}$ is an arbitrary bitangent line, then its complex conjugate is $b_{\sigma(i)\sigma(j)}$, and the set $I = \{i,\sigma(i),j,\sigma(j)\}$ corresponds to one of the three same index sets. This proves the claim by \Cref{thm:steiner}.
\end{proof}

With this knowledge we carry out Dixon's Algorithm with the following input data: Let $[D]$ be a strongly real 2-torsion point that cannot be represented as a difference of real bitangent lines. Let $Q$ be the associated quadratic representation. By choice, $\cC_Q(\R)$ does not contain any singular conics, since they would correspond to products of bitangent lines. Now let $a_{00}$ be an element $\cC_Q(\R)$. We make some specific choices in the algorithm before we later prove in \Cref{lem:indep} that the result is independent of that choice.

\begin{lem}
    The conic $a_{10}$ in Dixon's Algorithm can be chosen as a product of two real lines.
\end{lem}

\begin{proof}
    Recall that $A$ is the set of four points $\cV(a_{00}) \cap \cV(F)$. These points appear as two pairs of complex conjugate points, and for both pairs there is a unique real line connecting the two points. Choose $a_{10}$ as the product of these two lines.
\end{proof}

\begin{thm} \label{thm:posQDR}
    Let $a_{00}, a_{10}$ be as above. Then either the matrix $Q$ returned by Dixon's Algorithm or the matrix $-Q$ is a positive quadratic representation of $F$.
\end{thm}

\begin{proof}
    By continuity of eigenvalues it suffices to prove that $Q$ is positive semidefinite in a single point. We will chose this point to lie on $V(a_{10})$. Consider the matrix
    \begin{equation*}
        Q = \begin{pmatrix}
            a_{00} & a_{01} \\
            a_{01} & a_{11}
        \end{pmatrix}
    \end{equation*}
    as constructed. We will analyze the topological types of the real zero sets of $a_{00}$ and $a_{11}$. Both $a_{00}$ and $a_{11}$ are real contact conics and as elements of $\cC_Q$ they are identified with the points $(1,0)$ and $(0,1)$. Now any path
    \begin{equation*}
        \gamma\colon [0,1] \to \R^2
    \end{equation*}
    satisfying $\gamma(0) = (1,0), \gamma(1) = (0,1)$ corresponds to a continuous deformation of $a_{00}$ to $a_{11}$ in $\cC_Q(\R)$. By smoothness of the elements in $\cC_Q(\R)$, all their real zero sets have the same topological type, i.e., they are either all empty or consist of one oval. In particular, the zero sets of $a_{00}$ and $a_{11}$ have the same type. If their zero set is an oval, then its complement in the real projective plane has two connected components. One of these components is orientable and the other is not. The deformation preserves the sign that the quadratic form takes on each of these components. There are two cases that we need to distinguish:
    \begin{enumerate}
        \item $\cV(a_{00})$ has empty real part. Then so has $\cV(a_{11})$. We may then choose $P$ arbitrary on the real locus of $\cV(a_{10})$. Since positivity is preserved by $\gamma$ this shows $\det(Q(P)) = a_{00}(P)a_{11}(P) > 0$ and hence definiteness.

        \item $\cV(a_{00})$ consists of a single oval. Then so does $\cV(a_{11})$. Because the real zero set of $\cV(a_{10})$ is not entirely contained in the interior of any one of these ovals there is a point $P$ in the plane $\pp^2_\R(\R)$ that lies in the interior of neither oval, but at which $a_{10}$ vanishes. Again, since $\gamma$ preserves positivity, this shows $\det(Q(P)) = a_{00}(P)a_{11}(P) > 0$ and hence definiteness.\qedhere
    \end{enumerate}
\end{proof}

\begin{rem}
    In retrospect it is apparent that the second case in the proof of \Cref{thm:posQDR} will actually never occur since if $a_{00}$ had non-empty real zero locus, then $Q$ would not be a positive quadratic representation.
\end{rem}

\begin{lem} \label{lem:indep}
    The choice of $a_{10}$ in Dixon's algorithm does not affect whether the resulting quadratic representation $Q$ is positive.
\end{lem}

\begin{proof}
    First note that given $a_{00}$ and $a_{10}$ the conic $a_{11}$ is already uniquely determined. Indeed, if through Max Noether's fundamental theorem we obtained identities
    \begin{align*}
        a_{10}^2 = a_{00}a_{11} + pF \\
        a_{10}^2 = a_{00}a_{11}' + p'F,
    \end{align*}
    then necessarily
    \begin{equation*}
        a_{00}(a_{11} - a_{11}') = (p-p')F
    \end{equation*}
    with $p - p' \in \C$, which proves that either $F$ is reducible or both sides are identically zero. By assumption on $F$ only the latter can occur, so $a_{11} = a_{11}'$.

    Since $\ell(A) = 2$ by \Cref{lem:validchoice}, if we replace $a_{10}$ by another conic $a_{10}'$ we necessarily have a relation
    \begin{equation*}
        {a_{10}'} = \lambda a_{00} + \mu a_{10}
    \end{equation*}
    for some $\lambda, \mu \in \R$. Without loss of generality we assume $\lambda = 1$. This new choice gives rise to the unique new quadratic representation
    \begin{equation*}
        Q' = \begin{pmatrix}
            1 & 0 \\
            1 & \mu
        \end{pmatrix} Q \begin{pmatrix}
            1 & 1 \\
            0 & \mu
        \end{pmatrix}.
    \end{equation*}
    In particular, due to Sylvester's law of inertia, the signature of $A'$ is the same as that of $A$ in every point of evaluation, which proves independence from choice of $a_{10}$.
\end{proof}

\begin{rem}
    It is essential that we consider smooth ternary quartic curves. For example, it is known that the polynomial
    \begin{equation*}
        F = (x_2^2 - x_0^2)^2 + (x_1^2 - x_0^2)^2
    \end{equation*}
    does not have a positive quadratic representation. This was first established by \cite[\S 7]{quarez}, and could alternatively be proven using methods similar to those employed for the Robinson polynomial in \Cref{sec:Rob}.
\end{rem}

\subsection{Quartic Curves: A Constructive Approach}

While the previous section proves the existence of positive quadratic representations and establishes their number, in this section we describe a method of obtaining a positive quadratic representation of a smooth positive ternary quartic, if one already knows a real linear determinantal representation. Since by \Cref{prop:realsymmLDR} such representations always exist for quartic curves, this section provides another proof of the existence of positive quadratic representations for arbitrary smooth positive ternary quartic curves $C = \cV(F)$.

Let $M$ be a real linear determinantal representation of $F$. As defined in \Cref{prop:bitmatrix} let $O = (O_i)_{i=1}^8$ be the Cayley octad associated with $M$ (without loss of generality we assume that $O_{2i} = \overline{O_{2i-1}}$ for $i = 1, \dots, 4$) and let $\cB$ be the bitangent matrix.

\begin{lem} \label{lem:cayleybasechange}
    There is a change of basis $T \in Gl_4(\R)$ such that
    \begin{equation*}
        T\cdot(O_1^t, O_2^t, O_3^t, O_4^t) = \begin{pmatrix}
            1 & 1 & 0 & 0 \\
            i & -i & 0 & 0 \\
            0 & 0 & 1 & 1 \\
            0 & 0 & i & -i
        \end{pmatrix}.
    \end{equation*}
    After application of this base change we obtain the real linear determinantal representation
    \begin{equation*}
        M' = T^{-t}MT^{-1} = \begin{pmatrix}
            \frac{1}{2} b_{12} & 0 & l_{13} & l_{14} \\
            0 & \frac{1}{2} b_{12} & l_{23} & l_{24} \\
            &  & \frac{1}{2} b_{34} & 0 \\
            & & 0 & \frac{1}{2} b_{34}
     \end{pmatrix}
    \end{equation*}
    with the entries $l_{jk}$ being uniquely determined by the following identities:
\begin{align*}
    l_{13} &= \frac{1}{4} (b_{13} + b_{24} + (b_{14} + b_{23})) \\
    l_{24} &= -\frac{1}{4} (b_{13} + b_{24} - (b_{14} + b_{23})) \\
    l_{23} &= -\frac{i}{4} (b_{13} - b_{24} + (b_{14} - b_{23})) \\
    l_{14} &= -\frac{i}{4} (b_{13} - b_{24} - (b_{14} - b_{23}))
\end{align*}
\end{lem}

\begin{proof}
    By \cite[Proposition 6.3.2]{dolgachev} no four of the points $O_i$ are coplanar, so the desired base change
    \begin{equation*}
        T\cdot(O_1^t, O_2^t, O_3^t, O_4^t) = \begin{pmatrix}
            1 & 1 & 0 & 0 \\
            i & -i & 0 & 0 \\
            0 & 0 & 1 & 1 \\
            0 & 0 & i & -i
        \end{pmatrix}
    \end{equation*}
    exists over $\C$.
    If we apply complex conjugation on both sides and multiply on the right with the (column) permutation matrix $(12)(34) \in S_4$, then we obtain the identity
    \begin{equation*}
        \overline{T}\cdot(O_1^t, O_2^t, O_3^t, O_4^t) = \begin{pmatrix}
            1 & 1 & 0 & 0 \\
            i & -i & 0 & 0 \\
            0 & 0 & 1 & 1 \\
            0 & 0 & i & -i
        \end{pmatrix}.
    \end{equation*}
    Uniqueness of $T$ then proves that $T \in Gl_4(\R)$.

    For the description of $M' = (l_{ij})_{i,j}$ we first denote $O_i' = O_i T^t$. Then we obtain identities
    \begin{equation*}
        O_i' M' (O_j')^t = b_{ij}.
    \end{equation*}
    for all $i,j$. In particular, from $1 \leq i,j \leq 2$ we obtain identities
    \begin{align*}
        0 &= (1,i,0,0)M'(1,i,0,0)^t = l_{11} + 2i l_{12} - l_{22} \\
        b_{12} &= (1,i,0,0)M'(1,-i,0,0)^t = l_{11} + l_{22}.
    \end{align*}
    Since $M'$ is defined over $\R$ we conclude from the first identity that $l_{12} = 0$ and $l_{11} = l_{22}$. From the second we then infer that $l_{11} = l_{22} = \frac{1}{2}b_{12}$. A straightforward calculation in the same spirit now gives the desired description of $M'$.
\end{proof}

Next we need to determine a real geometrically irreducible contact conic. To find such conics we relate them to the system of contact cubics associated with $M$.

\begin{lem}
    Let $\Sigma(D)$ be a Steiner complex of type 1 with index set $I = \{i,j\}$, and let $Q$ be the quadratic representation associated with $D$. Then
    \begin{equation*}
        b_{ij}\cC_Q = \{ \lambda^t M^{\adj}\lambda \mid \lambda \in \langle O_i, O_j \rangle^\perp \}.
    \end{equation*}
    In particular, if $O_j = \overline{O_i}$, then $\cC_Q$ contains real contact conics, and each of them is geometrically integral.
\end{lem}

\begin{proof}
    See \cite[Lemma 6.7]{plaumann} for the statement that $\cC_Q$ contains real contact conics. From \Cref{lem:geomirr} we know that they are geometrically integral.
\end{proof}

Now consider the matrix $M'$ and pass to its adjugate. It is of the form
\begin{equation*}
    (M')^{\adj} = \left(\begin{array}{cc}
        \frac{1}{2}b_{34} \begin{pmatrix}
         \frac{1}{4} b_{12}b_{34} - l_{23}^2 - l_{24}^2 & l_{13}l_{23} + l_{14}l_{24} \\
         l_{13}l_{23} + l_{14}l_{24} & \frac{1}{4} b_{12}b_{34} - l_{13}^2 - l_{14}^2
        \end{pmatrix} & * \\
        * & *
    \end{array}\right).
\end{equation*}
From this matrix we not only obtain irreducible contact conics, but in fact the whole quadratic representation $Q$.

\begin{thm}
    The matrix
    \begin{equation*}
        Q = - \begin{pmatrix}
         \frac{1}{4} b_{12}b_{34} - l_{23}^2 - l_{24}^2 & l_{13}l_{23} + l_{14}l_{24} \\
         l_{13}l_{23} + l_{14}l_{24} & \frac{1}{4} b_{12}b_{34} - l_{13}^2 - l_{14}^2
        \end{pmatrix}
    \end{equation*}
    is a positive quadratic representation of $f$.
\end{thm}

\begin{rem}
    The bottom right corner of $(M')^{\adj}$ is the matrix
    \begin{equation*}
        \frac{1}{2}b_{12} \begin{pmatrix}
         \frac{1}{4} b_{12}b_{34} - l_{13}^2 - l_{23}^2 & l_{13}l_{14} + l_{23}l_{24} \\
         l_{13}l_{23} + l_{14}l_{24} & \frac{1}{4} b_{12}b_{34} - l_{14}^2 - l_{24}^2
        \end{pmatrix},
    \end{equation*}
    which similarly gives a positive quadratic representation.
\end{rem}

\begin{proof}
    By expansion of determinants it is verified that
    \begin{equation*}
        \det(Q) = \det\begin{pmatrix}
            \frac{1}{2} b_{12} & 0 & l_{13} & l_{14} \\
            0 & \frac{1}{2} b_{12} & l_{23} & l_{24} \\
            &  & \frac{1}{2} b_{34} & 0 \\
            & & 0 & \frac{1}{2} b_{34}
     \end{pmatrix} = F \geq 0.
    \end{equation*}
    Hence $Q$ is a quadratic representation of $F$. To see that it is positive it suffices to prove that all leading principal minors are positive in a real point where $f$ does not vanish (i.e. any real point). Since $F$ is positive by assumption it thus suffices to prove the inequality
    \begin{align*}
        l_{23}^2 + l_{24}^2 - \frac{1}{4} b_{12}b_{34} > 0
    \end{align*}
    in any real point of our choice.

    Using the identites from \Cref{lem:cayleybasechange} we may simplify the left hand side of the inequality to
    \begin{align*}
        4(l_{23}^2 + l_{24}^2 - \frac{1}{4} b_{12}b_{34}) &= (b_{13} - b_{23})(b_{24}-b_{14}) - b_{12}b_{34} = \\
        &= (b_{13} - b_{23})(\overline{b_{13}} - \overline{b_{23}}) - b_{12}b_{34} = \\
        &= ||b_{13} - b_{23}||^2 - b_{12}b_{34}
    \end{align*}
    with pointwise absolute value $|| \cdot ||$ taken in $\C$.
    
    The polynomial $b_{12}b_{34}$ vanishes in the union of two lines in $\pp^2_\R$, the other has a finite set of vanishing points. Thus the inequality is true for every point in $\cV(b_{12}b_{34})$ except maybe finitely many. This concludes the proof.
\end{proof}

\begin{ex}
Consider the linear determinantal representation of a positive quartic given in \cite[Example 5.12]{plaumann}:

{\scriptsize{
\begin{equation*}
    M = \begin{pmatrix}
        52 x_0 + 12x_1 - 60x_2 & -26x_0 - 6x_1 + 30x_2 & 48x_2 & 48x_1 \\
        -26x_0 - 6x_1 + 30x_2 & 26x_0 + 6x_1 - 30x_2 & -6x_0 + 6x_1 - 30x_2 & -45x_0 - 27x_1 -21x_2 \\
        48x_2 & -6x_0 + 6x_1 - 30x_2 & -96x_0 & 48x_0 \\
        48x_1 & -45x_0 - 27x_1 -21x_2 & 48x_0 & -48x_0
    \end{pmatrix}
\end{equation*}
}}

The associated Cayley octad is the set of points
\begin{equation*}
    O^t = \begin{pmatrix}
        i & -i & 0 & 0 & * & & & \\
        1+i & 1-i & 0 & 0 & & * & & \\
        0 & 0 & i & -i & & & * & \\
        0 & 0 & 1+i & 1-i & & & & *
    \end{pmatrix}
\end{equation*}
and the base change from \Cref{lem:cayleybasechange} is
\begin{equation*}
    T = \begin{pmatrix}
        -1 & 1 & 0 & 0 \\
        1 & 0 & 0 & 0 \\
        0 & 0 & -1 & 1 \\
        0 & 0 & 1 & 0
    \end{pmatrix}.
\end{equation*}
This lets us compute the positive quadratic representation
\begin{equation*}
    Q = \begin{pmatrix}
        q_{11} & q_{21} \\
        q_{21} & q_{22}
    \end{pmatrix}
\end{equation*}
with
\begin{align*}
    q_{11} &= 979x_0^2-726x_0x_1+195x_1^2+126x_0x_2-174x_1x_2+75x_2^2 \\
    q_{21} &= -771x_0^2+6x_0x_1+189x_1^2-774x_0x_2+198x_1x_2-99x_2^2 \\
    q_{22} &= 979x_0^2+810x_0x_1+195x_1^2+942x_0x_2+546x_1x_2+507x_2^2.
\end{align*}
\end{ex}

Finally, we argue that with the information of a real linear determinantal representation we may construct every positive quadratic representation of $F$. For the subsequent discussion we introduce the following standing assumptions, conventions and notations:
\begin{enumerate}
    \item In general we will work with divisors instead of line bundles on $C$, i.e. we identify line bundles $\cL$ with the unique divisor classes $[D]$ on $C$ satisfying $\cL \cong \cO(D)$ (see \cite[Corollary II.6.16]{Hart77}).

    \item By $H$ (as in \underline{h}yperplane) we will denote the divisor of zeros of a linear form on $C$.
    
    \item By $T_0$ (as in \underline{t}heta characteristic with \underline{$0$} global sections) we denote a representative of a strongly real even theta characteristic on $C$. Such a theta characteristic exists and $[T_0 + H]$ it defines a real linear determinantal representation $M$ of $F$ by \Cref{prop:realsymmLDR}. We assume $T_0 + H \sim T \geq 0$ with $T$ being invariant under complex conjugation.

    \item By $B$ (as in \underline{b}itangent line) we denote the unique effective representative of a strongly real odd theta characteristic on $C$. If further divisors of the form $B_i$ occur, they will have the same meanig (without necessarily the reality condition).

    \item By $A$ (as in \underline{a}ny other letter) we denote an effective and real representative of the strongly real divisor class $[T + B - H]$.
\end{enumerate}

In this list of assumptions we made some choices that are a priori unjustified. The next lemma proves that they are valid.

\begin{lem} \label{lem:validchoice}
    All choices in the above list are possible. Precisely we have the following:
    \begin{enumerate}
        \item All even theta characteristics on $C$ are non-vanishing. Hence, $T_0$ exists as desired.

        \item $\ell(T_0 + H) \geq 1$. Hence, $T$ exists as desired.

        \item $\ell(T + B - H) = 2$. Hence, $A$ exists as desired.
    \end{enumerate}
\end{lem}

\begin{proof}
    The first statement is a well-established fact for smooth quartic curves. For the second statement we use that
    \begin{equation*}
        \ell(T_0 + H) \geq \deg(T_0 + H) - g + 1 = 4.
    \end{equation*}
    Finally, we claim that
    \begin{equation*}
        \ell(A) = \deg(T + B - H) - g + 1 = 2,
    \end{equation*}
    which implies in particular that $A$ exists as desired. We prove this equality using Clifford's theorem. Assuming that $\ell(A) \geq 3$, then $2(\ell(A) - 1) \leq \deg(A) = 4$, so $\ell(A) = 3$ and equality holds. But this is only possible if $A$ is either canonical or the zero divisor, both of which is obviously not the case.
\end{proof}

\begin{rem}
    Given the linear determinantal representation $M$ the computation of $A$ can be done as outlined in  \cite[Lemma 6.7]{plaumann}). This uses that $2T$ is the intersection divisor of $C$ with a curve $C'$ given by a polynomial of the form $\lambda^t (M^{\adj}) \lambda$.
\end{rem}

The divisor class $[A - H]$ is a 2-torsion point on $J(C)$ and we claim that it induces a positive quadratic representation $Q$ of $F$.

\begin{rem}
    If $B_i, 1 \leq i \leq 4$ denote the contact divisors of the four real bitangent lines of $C$, and if $A_i \sim T_0 - B_i$ in an analogous way as above, then no two of these divisors are linearly equivalent. Indeed, if this were the case, then we had
    \begin{equation*}
        T_0 - B_i \sim T_0 - B_j \iff B_i \sim B_j \iff B_i - B_j \sim 0,
    \end{equation*}
    a contradiction. In particular, our setup lets us construct four different quadratic representations.
\end{rem}

\begin{lem} \label{lem:geomirr}
    Under our standing assumptions the system $\cC_Q(\R)$ does not contain any geometrically reducible conic sections. In particular, all elements of $\cC_Q(\R)$ are smooth conics.
\end{lem}

\begin{proof}
    The system $\cC_Q$ contains exactly six contact conics that are reducible over $\C$, namely the products of the six pairs of bitangent lines that lie in the Steiner complex $\Sigma(D)$ where $D = A -H \sim T_0 - B$. It suffices to prove that none of these six products is real. For that assume that contrarily
    \begin{equation*}
        D \sim B_1 - B_2
    \end{equation*}
    for two bitangent lines $B_1$ and $B_2$ that are either both real or related by complex conjugation. If $B_1, B_2$ is a pair of complex conjugate bitangent lines, then their difference is linearly equivalent to the difference of two real bitangent lines by \Cref{lem:real_azyg}. Thus without loss of generality we restrict to both $B_1, B_2$ being real bitangent lines.     
    We then write $D \overset{\text{Def}}{\sim} T + B - 2H \sim T_0 - B$ and observe that we get
    \begin{align*}
        T_0 - B \sim B_1 - B_2 &\iff T \sim T_0 + H \sim B + B_1 + B_2 \\
        &\iff T_0 \sim B + B_1 + B_2 - H.
    \end{align*}
    
    From \Cref{thm:steiner}, considering the Steiner complex $S^{(2)}_I$ with $I = \{1,2,3,4\}$, we infer that the contact points of the four real bitangent lines with $C$ lie on a plane conic. In particular, since $B, B_1, B_2$ are the contact divisors of three distinct such lines (note that if $B = B_1$ or $B = B_2$, then necessarily $T_0 \sim B_1$ or $T_0 \sim B_2$, a contradiction) and if $B'$ is the contact divisor of the last, we get
    \begin{equation*}
        B + B_1 + B_2 - H \sim B + B_1 + B_2 + B' + B' - 2H \sim B'.
    \end{equation*}
    Thus we get $T_0 \sim B'$, a contradiction to $\ell(T_0) = 0$.
\end{proof}

\begin{cor}
    The divisor $A$ gives rise to a positive quadratic representation of $F$.
\end{cor}

\begin{proof}
    We have shown in \Cref{lem:geomirr} that the hypotheses of \Cref{thm:posQDR} are satisfied. Thus the claim follows.
\end{proof}

\section{Concluding Remarks}
It is natural to ask whether our results for smooth quartic curves generalize to curves of higher degree. It is a quick observation that the set of nonnegative polynomials of degree $2d$ admitting a positive quadratic representation contains an open set in the Euclidean topology. To see this, one starts with the positive quadratic representation
\begin{equation*}
    F = (x_0^2 + x_1^2 + x_2^2)^d = \det((x_0^2 + x_1^2 + x_2^2) I_d)
\end{equation*}
and slightly perturbs this representation.
It is less apparent how large this set of polynomials is and we conclude with the following conjecture.

\begin{con}\label{con1}
    The set of nonnegative forms of degree $2d$ admitting a positive quadratic representation is dense (with respect to the Euclidean topology) in the set of all nonnegative forms of degree $2d$.
\end{con}

\appendix

\section{Blow-ups of isolated real nodes} \label{app:realblowup}
Let $R=\R[x,y]_{\fm}$ be the localization of $\R[x,y]$ at the maximal ideal $\fm=(x,y)$. Let
\begin{equation*}
    F=F_2+F_3+\ldots+F_n
\end{equation*}
with $F_i\in R[x,y]$ homogeneous of degree $i$ and $F_2\neq0$ irreducible. Note that this implies that $F$ itself is irreducible and hence $A=R/(F)$ is an integral domain. Further note that $A$ is a local ring with maximal ideal generated by $x$ and $y$ which, by abuse of notation, we also denote by $\fm$. We will explicitly describe the integral closure of $A$. To this end consider the ring extension $S=R[z]/(y-xz)$ of $R$ and
\begin{equation*}
    F'=F_2(1,z)+xF_3(1,z)+\ldots+x^{n-2}F_n(1,z)\in S.
\end{equation*}
We will prove that $A'=S/(F')$ is the integral closure of $A$. We note that since $F=x^2F'$, the natural ring homomorphism $R\to A'$ factors through
\begin{equation*}
    \pi\colon A\to A'.
\end{equation*}

We will use the following statements.

\begin{prop}\label{prop:blowup}
The following are true:
\begin{enumerate}
    \item The ring homomorphism $ \pi\colon A\to A'$ is injective.
    \item $A'$ is generated as an $A$-module by $z$.
    \item $ax\in A$ for every $a\in A'$.
\end{enumerate}
Note that in (2) and (3) we consider $A$ as a subring of $A'$ via $\pi$.
\end{prop}

\begin{proof}
 This is the statement of \cite[Proposition~1]{fultonblowup}. Note that in \cite{fultonblowup} they assume the ground field to be algebraically closed but the proof of \cite[Proposition~1]{fultonblowup} works over an arbitray field.
\end{proof}

\begin{cor}\label{cor:AinN}
$A'$ is contained in the normalization of $A$.
\end{cor}

\begin{proof}
 Because $A'$ is generated by $z$ as an $A$-module, there are $a_1,a_2\in A$ such that $z^2=a_1+a_2z$ which shows that $z$ and hence $A'$ is integral over $A$. Furthermore $A'$ is contained in the quotient field of $A$ because $y=xz$.
\end{proof}

Let $\fn$ be the ideal of $A'$ that is generated by $x$ and $y$. We have
\begin{equation*}
    A'/\fn=\R[z]/F_2(1,z)\cong\C
\end{equation*}
because $F_2(1,z)$ is irreducible of degree two. This implies that $A'$ is a local ring with maximal ideal $\fn$.

\begin{lem}\label{lem:Areg}
 $A'$ is a regular local ring.
\end{lem}

\begin{proof}
 The Jacobian matrix of the two equations $F'$ and $y-xz$ is
 \begin{equation*}
     \begin{pmatrix}
      -z& 1& -x\\
      \frac{\partial F'}{\partial x}& 0 & \frac{\partial F'}{\partial z}
     \end{pmatrix}.
 \end{equation*}
 It has rank $2$ modulo $\fn$ since $F_2(1,z)$ is irreducible and thus separable.
\end{proof}

\begin{cor}
 $A'$ is the normalization of $A$.
\end{cor}

\begin{proof}
By \Cref{cor:AinN} the ring $A'$ is contained in the normalization of $A$ and by \Cref{lem:Areg} $A'$ is regular and therefore integrally closed. This implies the claim.
\end{proof}

\begin{prop}\label{prop:Aexpl}
 We have
 \begin{equation*}
     A=\{f\in A'\mid f\equiv c\mod\fn\textrm{ for some }c\in\R\}.
 \end{equation*}
\end{prop}

\begin{proof}
Since $\fm=\fn\cap A$ and $A/\fm=\R$, we have the inclusion ``$\subseteq$''.
In order to prove the other inclusion, let $f$ an element of the right-hand side. By \Cref{prop:blowup} there are $a_1,a_2\in A$ such that $f=a_1+a_2z$. Reduced modulo $\fn$ this equals to
\begin{equation*}
    a_1(0,0)+a_2(0,0)\cdot \bar{z}.
\end{equation*}
Since $a_1(0,0)\in\R$ and $\bar{z}\not\in\R$, we must have $a_2(0,0)=0$, i.e. $a_2\in\fm$. Let $h_1,h_2\in A$ such that $a_2=xh_1+yh_2$. Then we have
\begin{equation*}
    a_2z=x\cdot(h_1z+h_2z^2)
\end{equation*}
which is an element of $A$ by \Cref{prop:blowup}. Thus $f\in A$ which proves the claim.
\end{proof}

\bibliographystyle{alpha}
\bibliography{Bibliography}

\section{Computational Methods for the Robinson Polynomial} \label{app:RobComp}

In this appendix we describe the computational methods employed in \Cref{sec:Rob}. All calculations were performed using \texttt{Macaulay2} (v.~1.24.05, \cite{M2}), and we give a detailed description of the algorithm. The main code file \texttt{Robinson} imports key functions from the file \texttt{RobinsonFunctions.m2} and comes preloaded with a classification of the representing types, stored in two files called \texttt{classes\_Robinson.txt} and \texttt{subclasses\_Robinson.txt}. It furthermore requires a text file \texttt{Parameters.txt}. Executing the main file then computes the quadratic determinantal representation for the 2-torsion point $a \in \mathbb{T}^{10}$ as defined in the latter file:
\begin{oframed}
\begin{verbatim}
-- Parameters.txt
-- Please define a as desired
a = class3#0;

------------------
///
available options:

Free Choice:
a in {1,i}^{10}
 -- admissible formats:
 -- a = (1, ..., 1)
 -- a = {1, ..., 1}

Preloaded Classes:

Det. Rep. of type M_0
a = class0#k -- 0..k..925

 -- Det. Rep. of type M_0^0
 -- sclass0#k -- 0..k..437

 -- Det. Rep. of type M_0^1
 -- sclass1#k -- 0..k..377

 -- Det. Rep. of type M_0^2
 -- sclass2#k -- 0..k..109

Det. Rep. of type M_1
a = class1#k -- 0..k..86

Det. Rep. of type M_2
a = class2#k -- 0..k..9

Det. Rep. of type M_3
a = class3#k -- k = 0
///
\end{verbatim}
\end{oframed}

\subsection{Initial Choices}

In terms of notation recall that $R$ denotes the Robinson polynomial, $X = \cV(R)$ and $Y$ is the singular locus of $X$. The morphism $\pi \colon Q \to X$ denotes the normalization of $X$. As a first step we compute the (unique) preimages of the elements of $Y$.

\begin{align*}
    q_1 &= \pi^{-1}(1:1:1) = [\frac{1+\sqrt{3}i}{2}:\frac{1-\sqrt{3}i}{2}:1] \\
    q_2 &= \pi^{-1}(1:1:-1) = [\frac{1-\sqrt{3}i}{2}:\frac{-1-\sqrt{3}i}{2}:1] \\
    q_3 &= \pi^{-1}(1:-1:1) = [\frac{-1+\sqrt{3}i}{2}:\frac{1+\sqrt{3}i}{2}:1] \\
    q_4 &= \pi^{-1}(1:-1:-1) = [\frac{-1+\sqrt{3}i}{2}:\frac{-1-\sqrt{3}i}{2}:1] \\
    q_5 &= \pi^{-1}(0:1:1) = [1:-1:\sqrt{2}i] \\
    q_6 &= \pi^{-1}(0:1:-1) = [1:1:\sqrt{2}i] \\
    q_7 &= \pi^{-1}(1:0:1) = [\sqrt{2}i:1:1] \\
    q_8 &= \pi^{-1}(1:0:-1) = [-\sqrt{2}i:-1:1] \\
    q_9 &= \pi^{-1}(1:1:0) = [1:\sqrt{2}i:1] \\
    q_{10} &= \pi^{-1}(1:-1:0) = [-1:-\sqrt{2}i:1]
\end{align*}

Choose the linear form $h = 2x_0 + x_1$ in $\pp^2_\R$. Its zero divisor on $X$ pulls back to the zero divisor of
\begin{equation*}
 h_\pi(x,y,z) = h((2x^2+z^2)z,x(x^2+2z^2),(x-z)(x+z)y)
\end{equation*}
cutting out the divisor $H_\pi = \cV(h_\pi)$ on $Q$, which satisfyies $H_\pi \cap \pi^{-1}(Y) = \emptyset$. Then the space of global sections of the sheaf $\cO_X(t)$ may be interpreted in $K(Q)$ as the $\R$-vectorspace spanned by
\begin{align*}
 \frac{x^iy^i}{h_\pi^t}, &\qquad i+j=3t \\
 \frac{x^iy^iz}{h_\pi^t}, &\qquad i+j=3t-1.
\end{align*}

\subsection{Setup}

Unless already installed, our code requires the package \texttt{RealRoots}. Furthermore, some of the functions employed are imported from a second file called \texttt{RobinsonFunctions.m2}, thus the code starts as follows:

\begin{oframed}
\begin{verbatim}
-- Computation of QDRs of the Robinson Polynomial
restart

needsPackage "RealRoots";
load "RobinsonFunctions.m2";
\end{verbatim}
\end{oframed}

In line with the code below we will denote the Robinson form by $Rob$ (to free up the letter $R$ for a ring). Since all points $q_i$ lie in the complement of $\cV(z)$ we will exclusively work in the affine chart $z \neq 0$, and we introduce the affine coordinates $v = \frac{x}{z}, w = \frac{y}{z}$. Furthermore we define the square roots $s2 = \sqrt{2}$ and $s3 = \sqrt{3}$ via the minimal polynomial $\mu^4 - 10\mu^2 + 1$ of $\sqrt{2} + \sqrt{3}$ over $\Q$.

To force Macaulay to interpret $R$ as a field we use the command \texttt{toField}, and we add the imaginary unit $i$ for later calculations. For comparison of coefficients we do not want Macaulay to interpret $i$ as a field element, hence the late adjunction.

\begin{oframed}
\begin{verbatim}
R = QQ[v,w,mu]/(ideal(mu^4 - 10*mu^2 + 1, v^2 + w^2 + 1));
    -- affine coordinate ring of Q with mu = sqrt(2) + sqrt(3)
    s2 = (mu^3-9*mu)/2; -- s2 = sqrt(2)
    s3 = (11*mu-mu^3)/2; -- s3 = sqrt(3)

F = frac R;
    -- F = K(Q);

K = toField(F[i]/(ideal(i^2+1)));
    -- i = sqrt(-1);
\end{verbatim}
\end{oframed}

Next we define the Robinnson form $Rob$ in the projective plane with coordinates $z_0, z_1, z_2$. As before we adjoin the element $\sqrt{2} + \sqrt{3}$, denoted by $\mu_z$.

\begin{oframed}
\begin{verbatim}
-- Robinson form defined in projective plane with coordinates zi
S = QQ[z0,z1,z2,muz]/(ideal(muz^4 - 10*muz^2 + 1));
Rob = z0^6 + z1^6 + z2^6 - z0^4*z1^2 - z0^4*z2^2 -z0^2*z1^4 -
      z0^2*z2^4 - z1^4*z2^2 - z1^2*z2^4 + 3*z0^2*z1^2*z2^2;
\end{verbatim}
\end{oframed}

Finally, once we have found a quadratic representation, we need to evaluate it in a point, which will produce a matrix over the field $\Q(\sqrt{2} + \sqrt{3})$, which we denote by $S2$.

\begin{oframed}
\begin{verbatim}
-- ring for algebraic psd condition
S2 = QQ[u];
MinPol = u^4-10*u^2+1;
E = subsets(3);
\end{verbatim}
\end{oframed}

Initially we work in the field $K = \Q(\sqrt{2}, \sqrt{3}, i)$ and define the normalization map. For later use we store the mapping data $z_i = X_i$ in lists $x$ and $z$. We also define the list $Q$ of preimages of the singular points of $Rob$, as well as the polynomial defining the pullback of a line not containing any of them.

\begin{oframed}
\begin{verbatim}
use K;

-- normalization map defined by zi = Xi
X0 = 2*v^2+1;
X1 = v^3+2*v;
X2 = (v^2-1)*w;

x = {X0,X1,X2};
z = {z0,z1,z2};

-- points of evaluation
q1 = {(1+s3*i)/2,(1-s3*i)/2};
q2 = {(1-s3*i)/2,(-1-s3*i)/2};
q3 = {(-1+s3*i)/2,(1+s3*i)/2};
q4 = {(-1+s3*i)/2,(-1-s3*i)/2};
q5 = {-s2*i/2,s2*i/2};
q6 = {-s2*i/2,-s2*i/2};
q7 = {s2*i,1}; 
q8 = {-s2*i,-1};
q9 = {1,s2*i};
q10 = {-1,-s2*i};

Q = {q1,q2,q3,q4,q5,q6,q7,q8,q9,q10};

-- divisor H, pullback of line section
H = 2*X0 + X1;
\end{verbatim}
\end{oframed}

\subsection{Calculations}

For explicit calculations we start by describing a method for finding a basis of the vector spaces
\begin{equation*}
    \{ f \in H^0(\cO_Q(tH)) \mid f(q_k) \in \R a_k \}
\end{equation*}
for a given tuple $a \in (\Z/2\Z)^{10}$.

First recall that an element $f \in H^0(\cO_Q(tH))$ is of the form
\begin{equation*}
    f(v,w) = \sum_{3t-1 \leq i+j \leq 3t} c_{ij} \frac{v^iw^j}{H^t}.
\end{equation*}

We first write down a basis for $H^0(\cO_Q(tH))$ (for $t = 0, \dots, 4$) as follows with the function \texttt{GG} (as in \underline{G}lobal \underline{G}enerators). Since we barely ever use the denominators of these basis vectors we simply exclude them from the output.

\begin{oframed}
\begin{verbatim}
-- global sections of O_P^2(t), i.e. of O_Q(tH)
GG = (H,t) -> (
    -- H ... polynomial defining divisor H
    -- t ... desired twist
    -- return: 1) basis of polynomials mod (v^2 + w^2 + 1)
    --            of degree 3t
    --         2) twist t

    k := 2*3*t+1; -- number of polynomials returned
    GlobGens := new MutableList from (k:0);
    
    for j in (0..(3*t)) do (
        GlobGens#j = v^(3*t-j)*w^j;
    );
    for j in (0..(3*t-1)) do (
        GlobGens#(j+3*t+1) = v^(3*t-j-1)*w^j;
    );
    
    return (toList GlobGens, t);
)
\end{verbatim}
\end{oframed}

The additional condition $f(q_k) = \lambda_k a_k$ for $k = 0, \dots 9$ defines a homogeneous system of 10 linear equations in the variables $(c_{ij}, \lambda_k)$ over $\C$. The number of equations doubles if we consider the system over $\R$ (splitting the original system in real and imaginary part). We build the corresponding matrix of the system in two steps. First we calculate $f(q_k)$ for all $k$ and split it into real and imaginary part. The result is a matrix of size $20 \times 10$. It is defined over $F$, so in particular we forget the existence of $i$.

For example, if we denote the above basis of $\cO_Q(H)$ by $(f_{ij})_{2 \leq i + j \leq 3}$ the matrix is defined over $\C$ as
\begin{equation*}
    (f_{ij}(q_k))_{i,j,k} = \begin{pmatrix}
        f_{30}(q_0) & f_{21}(q_0) & \dots & f_{02}(q_0) \\
        f_{30}(q_1) & f_{21}(q_1) & \dots & f_{02}(q_1) \\
        \vdots & \vdots & \ddots & \vdots \\
        f_{30}(q_9) & f_{21}(q_9) & \dots & f_{02}(q_9) \\
    \end{pmatrix}.
\end{equation*}
Over $\R$ this becomes
\begin{equation*}
    \binom{\text{Re}(f_{ij}(q_k))}{\text{Im}(f_{ij}(q_k))}_{i,j,k} = \begin{pmatrix}
        \text{Re}(f_{30}(q_0)) & \text{Re}(f_{21}(q_0)) & \dots & \text{Re}(f_{02}(q_0)) \\
        \text{Im}(f_{30}(q_0)) & \text{Im}(f_{21}(q_0)) & \dots & \text{Im}(f_{02}(q_0)) \\
        \vdots & \vdots & \ddots & \vdots \\
        \text{Re}(f_{30}(q_9)) & \text{Re}(f_{21}(q_9)) & \dots & \text{Re}(f_{02}(q_9)) \\
        \text{Im}(f_{30}(q_9)) & \text{Im}(f_{21}(q_9)) & \dots & \text{Im}(f_{02}(q_9)) \\
    \end{pmatrix}.
\end{equation*}

To introduce the condition imposed by $a \in (\Z/2\/)^{10}$ we have to extend this ``evaluation matrix'' to the right by columns of the form $- e_{2k}$ (if $f(q_k) \in \R$) or $- e_{2k+1}$ (if $f(q_k) \in \R i$). The function \texttt{aGens} (as in \underline{Gen}erator\underline{s} of $\cL_{\underline{a}}$) defines the matrix and calculates its kernel, i.e. a basis of $H^0(X, \cO_X(t))$. Since for later calculations the denominators of the base vectors are irrelevant, the function will not return them.

\begin{oframed}
\begin{verbatim}
-- Generators of L_a(t)
aGens = (GlobGens,Q,a) -> (
    -- GlobGens ... given set of rational functions f1, f2, ...
    --              = return pair of GG
    -- Q ... given set of points for evaluation
    -- a ... input vector of length 10 with entries in 1,i
    -- return: enumerators of generators f of O_Q(tH)
    --         satisfying f(qi)*1/H^t = RR*ai

    t = GlobGens#1;
    l = length(GlobGens#0);
    
    EvMat := mutableMatrix(F,20,l + 10);
    -- matrix of zeros with size of return
    
    for k in (0..length(Q)-1) do (
        ev := (GlobGens#0) / (section ->
            sub(section, {v=>(Q#k)#0, w=>(Q#k)#1})*
            sub(1/(H^t), {v=>(Q#k)#0, w=>(Q#k)#1}));

            -- real part of evaluation
            evR := ev / (section -> sub (section,{i=>0}));
            -- imaginary part of evaluation
            evI := (ev-evR) / (section -> sub(section,{i=>1}));
            
        for j in (0..l-1) do (
            EvMat_(2*k,j) = evR#j;
            EvMat_(2*k+1,j) = evI#j;

        if a#k == 1 then EvMat_(2*k,l + k) = -1
        else EvMat_(2*k + 1,l + k) = -1; -- else if a#k == i
        );
    );

    aGenerators := matrix{GlobGens#0} *
                   (gens(kernel(matrix(EvMat))))^{0..6*t};
    
    return aGenerators;
)
\end{verbatim}
\end{oframed}

For later use we make the following one-time calculation:
\begin{oframed}
\begin{verbatim}
-- generators of cO_X(1)
ONE = (10:1);
G = aGens(GG(H,1), Q, ONE);
\end{verbatim}
\end{oframed}

We now proceed by definnig the function that for given $a \in (\Z/2\Z)^{10}$ will return a corresponding real determinantal representation. Given $a$ we can make a case distinction based on the dimension of $h^0(X, \cL_a(1))$. Depending on the question we want to answer we can later preselect the tuples $a$ of interest.

\begin{oframed}
\begin{verbatim}
-- case distinction based on h^0(C, La(1))
cases = a -> (
    -- a ... input vector of length ten with entries in 1,i
    -- return: 1) H^0(X, La(1));
    --         2) h^0(X, La(1));
    
    Gens1 := aGens(GG(H,1),Q,a);
    return (Gens1, numColumns(Gens1));
)
\end{verbatim}
\end{oframed}

The function \texttt{LaComputation} (as in \underline{computation} of the determinantal representation associated with $\underline{\cL_a}$) first computes \texttt{aGens}(\texttt{GG}($H$, $1$), $Q$, $a$) and performs the case distinction. Furthermore, and initially just for $t = 0, 1$ we store the calculated generators \texttt{aGens}(\texttt{GG}($H$,$t$),$Q$,$a$).

\begin{oframed}
\begin{verbatim}
-- computation of the det. rep. associated with La
LaComputation = a -> (
    print a;
    c := cases(a);

    aGenList := {aGens(GG(H,0),Q,a), c#0};

    if c#1 == 3 then (
        FinalGens := aGenList#0;
        return DetRep(FinalGens, c#1);
    )

    ...
\end{verbatim}
\end{oframed}

The first case is when $h^0(X, \cL_a(1)) = 3$. Then it is not hard to a priori establish that the only possible matrix representation is of the representing type $M_3$ In this case the function \texttt{DetRep} calculates $M$ as follows.

\begin{oframed}
\begin{verbatim}
-- calculation of det. rep.
DetRep = (FinalGens, case) -> (
    -- FinalGens ... suitable list of generators
    -- case ... number of global generators of La(1) 
    -- return: 1) det. rep. of Rob
    --         2) degree type of det. rep.

    -- Trivial Case (i.e. there is a global section of F(0))
    if case == 3 then (
        Ansatz := sub(FinalGens * PowerX(x,6),R);
        -- matrix entry of degree 6 in x
        
        LinEq := sub((coefficients(Ansatz, Variables => {s,t},
                                   Monomials => Mon(s,t,6)))#1,F);
        
        return (PowerX(z,6) * sub(gens kernel LinEq,S),
                matrix{{6}});
    )

    ...
\end{verbatim}
\end{oframed}

The \texttt{Ansatz} writes the matrix $M$ as a $1 \times 1$ matrix containing a degree 6 polynomial with variables as coefficients. The coefficients of the polynomial in $M$ are then calculated by a comparison of coefficients.

The next case is when $h^0(X, \cL_a(1)) = 2$, meaning that the representing type will be $M_2$. This, however, poses the following problem: In theory we would like to determine the matrix $M$ from a short exact sequence
\begin{equation*}
    0 \to \cO(-4)^{\oplus2} \oplus \cO(-1) \overset{M}{\to} \cO^{\oplus2} \oplus \cO(-3) \overset{N}{\to} \cL_a(1) \to 0
\end{equation*}
by postulating $NM = 0$. To do so, however, we would need at least three global sections of $\cL_a(1)$, which we do not have. Instead we have to take higher twists until we have the sequence
\begin{equation*}
    0 \to \cO(-1)^{\oplus2} \oplus \cO(2) \overset{M}{\to} \cO(3)^{\oplus2} \oplus \cO \overset{N}{\to} \cL_a(4) \to 0.
\end{equation*}
Clearly we now have sufficiently many global sections, and we have natural global generators, since both $\cO(3)$ and $\cO$ are globally generated. We now have to choose suitable such generators in order to determine $N$. We already know two of them (from $\cL_a(1)$), and the last must be the one arising from the second part of the direct sum.

In practice we will do the following. We know that $H^0(X, \cL_a(1)) = \R f_1 + \R f_2$ for some $f_1, f_2$. We then consider the injective map
\begin{equation*}
    0 \to \langle f_i g | g \in G, i = 1,2 \rangle \to H^0(X, \cL_a(2)).
\end{equation*}
Recall that $G$ is a basis of $H^0(X, \cO_X(1))$. As new potential generators we add a basis of the cokernel $K_2$ of this map. However, $K_2 = 0$, because otherwise we had an entry of degree 3 in $M$. Ultimately we have to consider
\begin{equation*}
    0 \to \langle f_i g | g \in G^3, i = 1,2 \rangle \to H^0(X, \cL_a(4)) \to K_4 \to 0
\end{equation*}
to find new generators in $K_4 \neq 0$. These correspond to the entries of degree 1 in $M$.

\begin{ex}
If $h^0(X, \cL_a(1)) = 1$ a similar problem arises. Consider a determinantal representation $M$ of the form
\begin{equation*}
    \begin{pmatrix}
        4 & 3 \\
        3 & 2
    \end{pmatrix}.
\end{equation*}
We have $H^0(X, \cL_a(1)) = \R f$ for some $f$, but we need at least two global sections to define $N$ in the sequence
\begin{equation*}
    0 \to \cO(-5) \oplus \cO(-4) \overset{M}{\to} \cO(-1) \oplus \cO(-2) \overset{N}{\to} \cL_a.
\end{equation*}
Thus we twist twice and obtain
\begin{equation*}
    0 \to \cO(-3) \oplus \cO(-2) \overset{M}{\to} \cO(1) \oplus \cO(0) \overset{N}{\to} \cL_a(2),
\end{equation*}
from which we may now compute $N$.

Note, however, that alone from $h^0(X, \cL_a(1)) = 1$ we cannot unambiguously decide that the representing type of the Robinson polynomial is $M_1^0$. It could be $M_1^1$. In this case the initial short exact sequence is
\begin{equation*}
    0 \to \cO(-5) \oplus \cO(-4) \oplus \cO(-3) \overset{M}{\to} \cO(-1) \oplus \cO(-2) \oplus \cO(-3) \overset{N}{\to} \cL_a,
\end{equation*}
and the second twist is
\begin{equation*}
    0 \to \cO(-3) \oplus \cO(-2) \oplus \cO(-1) \overset{M}{\to} \cO(1) \oplus \cO(0) \oplus \cO(-1) \overset{N}{\to} \cL_a(2).
\end{equation*}
Critically, since $\cO(-1)$ is not globally generated, there is no natural way of globally generating $\cL_a$ through $\cO(1) \oplus \cO(0) \oplus \cO(-1)$, and we have to compute the third twist as well. Then, finally, from
\begin{equation*}
    0 \to \cO(-2) \oplus \cO(-1) \oplus \cO \overset{M}{\to} \cO(2) \oplus \cO(1) \oplus \cO \overset{N}{\to} \cL_a(3) \to 0
\end{equation*}
we obtain sufficiently many global generators.
\end{ex}

\begin{ex}
Consider a determinantal representation $M$ of the form $M_0^k$. We have $h^0(X, \cL_a(2)) = 3$, so there would be no problem in assuming that $M$ is of the form $M_3^0$. However, the three generators obtained in $\cL_a(2)$ need not globally generate $\cL_a(2)$ (and might even have common base points). In that case we also have to twist the initial sequence
\begin{equation*}
    0 \to \cO(-4)^{\oplus3} \oplus \cO(-3)^{\oplus k} \to \cO(-2)^{\oplus3} \oplus \cO(-3)^{\oplus k} \to \cL_a.
\end{equation*}
a third time, which results in the surjection
\begin{equation*}
    \cO(1)^{\oplus3} \oplus \cO^{\oplus k} \to \cL_a(3),
\end{equation*}
which now lets us find suitable global generators.
\end{ex}

\begin{rem}
    In practice we will list all base point free choices of global generators and test, which representations yield the Robinson polynomial. If this does not suffice (because the chosen generators do not globally generate $\cL_a(t)$ for suitable $t$), then we equally arbitrarily add another generator. This procedure ensures that we do not end up with representations containing superfluous data, e.g. representations of the form
    \begin{equation*}
        M = \begin{pmatrix}
            M' & 0 \\
            0 & 1
        \end{pmatrix},
    \end{equation*}
    which clearly admit a smaller representation.
\end{rem}

The relevant code to realize the embeddings
\begin{equation*}
    0 \to \langle f_i g | g \in G, i = 1,2 \rangle \to H^0(X, \cL_a(t+1)).
\end{equation*}
for $f_i \in \cL_a(t)$ is the following. The functions are called \texttt{EMB} (as in \underline{Emb}ed) and \texttt{EMBQuot} (as in \underline{Emb}ed and take a \underline{Quot}ient).

\begin{oframed}
\begin{verbatim}
-- Embed matrix of generators into next higher twist
EMB = Gens -> (
    -- Gens ... output pair of GG
    --          Gens0 = g1, g2, g3, ...
    -- return: g1*G1, g1*G2, g1*G3, g2*G1, g2*G2, g2*G3, ...
    --         in form of a row matrix
    
    Gens0 := (entries(Gens))#0; -- entries of matrix as list
    EMBED := Gens0 / (g -> sub(g*G,R));
    
    M := Zeros(1,0);
    for k in (0..length(Gens0)-1) do (
        M = M|(EMBED#k);
    );
    return M;
)
\end{verbatim}
\end{oframed}

This is followed by taking a quotient to remove redundant information. As the first input we take the list of bases of $H^0(X, \cL_a(t'))$ for $t' = 0, \dots, t$ with $t$ minimal such that this space is non-trivial. We want to apply our embedding procedure to the last entry of that list, and the argument \texttt{twist} determines how far we shift the twist. This procedure is iterative.

We begin by applying the function \texttt{EMB} to the generators of $H^0(X, \cL_a(t)$. This produces a subspace \texttt{SubMod1} of $H^0(X, \cL_a(t+1)) =$ \texttt{Mod1}. We interpret these spaces as free $R$-modules (to avoid computations with fractions), and compute a basis of the quotient Quot1.

If needed we iterate to pass to the next higher twist, using the elements of the first embedding and also those of $H^0(X, \cL_a(t+1))$ that are not yet covered by them (and so on).

\begin{oframed}
\begin{verbatim}
-- quotient of embedding H^0(X, La(t)) in H^0(X, La(t + twist))
EMBQuot = (aGenList,t,twist) -> (
    -- aGenList ... list of generators in different twists
    -- t ... lowest t such that h^0(X, La(t)) > 0
    -- twist ... target shift
    
    if twist == 0 then return aGenList

    -- if twist >= 0
    else (
    EMBED1 := EMB(aGenList#t);
    
    SubMod1 := sub((coefficients(EMBED1,Variables => {v,w},
                                 Monomials => Mon(v,w,t+1)))#1,R);
    Mod1 := sub((coefficients(sub(aGenList#(t+1),R),
                              Variables => {v,w},
                              Monomials => Mon(v,w,t+1)))#1,R);
    Quot1 := matrix{Mon(v,w,t+1)}*
                    mingens((image Mod1)/(image SubMod1));

    ------------------------------------------------------------
    if twist >= 2 then (
    
    EMBED2 := EMB(EMBED1 | Quot1);
        -- morally equal to EMB(EMB(aGenList#t) | aGenList#(d+1));
        
    SubMod2 := sub((coefficients(EMBED2,Variables => {v,w},
                                 Monomials => Mon(v,w,t+2)))#1,R);
    Mod2 := sub((coefficients(sub(aGenList#(t+2),R),
                              Variables => {v,w},
                              Monomials => Mon(v,w,t+2)))#1,R);
    Quot2 := matrix{Mon(v,w,t+2)}*
                    mingens((image Mod2)/(image SubMod2));
    
    if twist == 3 then (
    
    EMBED3 := EMB(EMBED2 | Quot2);
    
    SubMod3 := sub((coefficients(EMBED3,Variables => {v,w},
                                 Monomials => Mon(v,w,t+3)))#1,R);
    Mod3 := sub((coefficients(sub(aGenList#(t+3),R),
                              Variables => {v,w},
                              Monomials => Mon(v,w,t+3)))#1,R);
    Quot3 := matrix{Mon(v,w,t+3)}*
                    mingens((image Mod3)/(image SubMod3));
    return(Quot1,Quot2,Quot3);
    )
    
    else return(Quot1,Quot2);
    )
    ------------------------------------------------------------
    
    else return({Quot1});
    );
)
\end{verbatim}
\end{oframed}

Now assume that we have determined our set of ``final generators''. The calculation of $M$ is then in general performed as follows. Given our set \texttt{FinalGens} of $N$ many final generators in different twists $\cL_a(d)$ the determinantal representation must be the kernel of a surjection
\begin{equation*}
    0 \to \bigoplus_{i=1}^N \cO_X(-t_i' + t + \text{ twist}) \overset{M}{\to} \bigoplus_{i=1}^N \cO_X(-t_i + t + \text{ twist}) \overset{\text{FinalGens}}{\to} \cL_a(t + \text{ twist}) \to 0.
\end{equation*}
Since $\cL_a$ is supported on $\det(M)$ and $GM$ must vanish on $C = \cV(Rob)$ the determinant $\det(M)$ is a multiple of $Rob$. This describes a system of linear equations on the coefficients of the entries of $M$ which can then be solved.

\begin{rem}
    The mentioned system might return more than one solution, which is due to the fact that it does not inherently detect whether the constructed matrix $M$ is invertible (for example the zero matrix would be an admissible solution based on the condition $GM = 0$).
\end{rem}

Having already dealt with the case that $h^0(X, \cL_a(1)) = 3$, the next case is $h^0(X, \cL_a(1)) = 2$, and a determinantal representation is calculated by the following code.

\begin{oframed}
\begin{verbatim}
LaComputation = a -> (
    
    ...
    
    else if c#1 == 2 then (
        aGenList = aGenList | {aGens(GG(H,2),Q,a),
                               aGens(GG(H,3),Q,a),
                               aGens(GG(H,4),Q,a)};

        FinalGens0 := aGenList#1;
        EMBED := EMBQuot(aGenList,1,3);
        for k in (1..3) do (
            FinalGens0 = FinalGens0 | EMBED#(k-1);
        );

        l1 = numColumns FinalGens0;

        evGens := ((entries FinalGens0)#0) /
                  (section -> transpose matrix{
                      Q / (q -> tv(sub(section,
                                   {v => q#0, w => q#1}) != 0))
                  });
        
        Choice := ChooseGens(evGens,{2, l1, l1}, {2,1,0});
        for ch in Choice do (
            FinalGens := FinalGens0_ch;
            DR := DetRep(FinalGens,c#1);
            if (numColumns (DR#0) > 1) then (
                return DR;
            )
            else null;
        );
    )

    ...
\end{verbatim}
\end{oframed}

The code of \texttt{LaComputation} for this case first calculates \texttt{EMBQuot}(\texttt{aGenList},1,3), and collects all (linearly independent) generators in a row matrix \texttt{FinalGens0}. It assigns to each of them a column matrix with entries in 0 and 1, depending whether it vanishes in $q_k$ or not, which allows to decide which of the generators have common base points in $Q$. From all these generators the function \texttt{ChooseGens} (as in \underline{Choose} \underline{Gen}erator\underline{s}) produces a list of all possible choices of generators (in this case fixing the first 2, corresponding to the entries of degree 4 in $M_2$) such that they have no base points within $Q$. From these choices of generators a determinantal representation is calculated, again via comparison of coefficients.

\begin{rem}
    The choice of generators is a heuristically motivated method. There is no guarantee a priori that this in fact returns a determinantal representation of $Rob$. When no determinantal representation is retreived we add further generators to our potential list of global generators, and try again until we succeed.
\end{rem}

Both remaining cases, i.e. when $h^0(X, \cL_a(1)) \leq 1$ work in the same way, and we refrain from an explicit walk-through.

Ultimately we may verify, which of the representations of type $M_0^0$ are positive. Using a corollary of Hermite's Theorem (\cite[Theorem 4.13]{BasuPollackRoyAlg}) we may do so symbolically.
\begin{Def}
    Let $p \in \R[x]$ be a monic polynomial of degree $n$ with roots $(\alpha_i)_i$ in $\C$ and let $q \in \R[x]$ be arbitraty. We call the number
    \begin{equation*}
        \nu_k(p,q) = \sum_i \alpha_i^kq(\alpha_i) \in \R
    \end{equation*}
    the $k$-th generalized Newton sum of $p$. The matrix
    \begin{equation*}
        \mathscr{H}(p,q) = (\nu_{i+j}(p,q))_{0 \leq i,j \leq n-1}
    \end{equation*}
    is called the generalized Hermite matrix of $p$ and $q$.
\end{Def}

\begin{thm}
    Let $p \in \R[x]$ be monic and let $q, q_1, \dots, q_m \in \R[x]$. Furthermore for $e \in {1,2}^m$ define
    \begin{equation*}
        q^e = \prod_{i=1}^m q_i^{e_i}.
    \end{equation*}
    Then we have
    \begin{align*}
        \operatorname{sign} \mathscr{H}(p,q) &= \sum_{\alpha \in \R, p(\alpha) = 0} \operatorname{sign} q(\alpha) \\
        |\{ \alpha \in \R \mid p(\alpha) = 0, q_1(\alpha) > 0, \dots, q_m(\alpha) > 0 \}| &= \frac{1}{2^m} \sum_e \operatorname{sign} \mathscr{H}(p,q^e).
    \end{align*}
\end{thm}

\begin{proof}
    See \cite[Theorem 1.3.33 and Proposition 1.3.36]{Scheiderer2024RAG}.
\end{proof}

In our setting the polynomials $q_1,q_2,q_3$ will be the leading principal minors of a given quadratic representation, and if they are all positive, then the representation is positive. But positivity of the minors depends on the exact value of $\sqrt{2} + \sqrt{3}$, which we only know by its minimal polynomial, which will take the role of $p$.

To apply this for a given quadratic representation of the Robinson polynomial we first compute the evaluation in one of its points, and we map the resulting matrix into the ring $S2 = \Q(\sqrt{2} + \sqrt{3})$. Then we proceed as follows: We use that all real roots of the minimal polynomial of $\sqrt{2} + \sqrt{3}$ lie between $-4$ and $4$. To compute the signatures of the Hermite matrices we use the implementation \cite{RealRootsSource} of Sturm's algorithm in \texttt{Macaulay2} as described in \cite[Theorem 3]{M2RealRoots}.

Finally, given a quadratic representation we may evaluate it in a point and the following function checks, if the resulting matrix is positive semidefinite.

\begin{oframed}
\begin{verbatim}
isPSD = M -> (
    Minor0 := det M_{0}^{0};
    Minor1 := det M_{0,1}^{0,1};
    Minor2 := det M_{0,1,2}^{0,1,2};
    sylv := {};
    if Minor0*Minor1*Minor2 == 0 then return false
    else (
        for sets in E do (
            e := {0,1,2} / (i -> tv(isSubset({i},sets))+1);
            qe := Minor0^(e#0)*Minor1^(e#1)*Minor2^(e#2);
            sylv = sylv | {SylvesterCount(MinPol, qe, -4, 4)};
        );
        if sum sylv == 0 then return false
        else return true;
    )
)
\end{verbatim}
\end{oframed}

\end{document}

%% file: Operators.tex
\usepackage{mathrsfs}

\newcommand{\Pic}{\operatorname{Pic}}

\newcommand{\Proj}{\operatorname{Proj}}

\newcommand{\Div}{\operatorname{Div}}

\newcommand{\adj}{\operatorname{adj}} 


%% file: Fonts.tex
\newcommand{\cB}{{\mathcal B}}
\newcommand{\cC}{{\mathcal C}}

\newcommand{\cF}{{\mathcal F}}

\newcommand{\cL}{{\mathcal L}}

\newcommand{\cO}{{\mathcal O}}

\newcommand{\cU}{{\mathcal U}}
\newcommand{\cV}{{\mathcal V}}

\newcommand{\fm}{{\mathfrak m}}
\newcommand{\fn}{{\mathfrak n}}

\newcommand{\C}{{\mathbb C}}
\newcommand{\R}{{\mathbb R}}
\newcommand{\pp}{\mathbb{P}}
\newcommand{\Q}{{\mathbb Q}}
\newcommand{\N}{{\mathbb N}}
\newcommand{\Z}{{\mathbb Z}}




%% file: QuartCurves.bbl
\def\cprime{$'$}
\begin{thebibliography}{GMSY24}

\bibitem[Bea00]{beauville}
Arnaud Beauville.
\newblock Determinantal hypersurfaces.
\newblock {\em Michigan Math. J.}, 48:39--64, 2000.
\newblock Dedicated to William Fulton on the occasion of his 60th birthday.

\bibitem[BPR06]{BasuPollackRoyAlg}
Saugata Basu, Richare Pollack, and Marie-Fran\c{s}oise Roy.
\newblock {\em Algorithms in Real Algebraic Geometry}.
\newblock Number~10 in Algorithms and Computation in Mathematics. Springer,
  2006.

\bibitem[B{\v{S}}20]{buckleysivic}
Anita Buckley and Klemen {\v{S}}ivic.
\newblock New examples of extremal positive linear maps.
\newblock {\em Linear Algebra Appl.}, 598:110--144, 2020.

\bibitem[Dix02]{Dixon1902LDR}
Alfred~Cardew Dixon.
\newblock Note on the reduction of a ternary quantic to a symmetrical
  determinant.
\newblock {\em Mathemtical Proceedings of the Cambridge Philosophical Society},
  11:350--351, 1902.

\bibitem[Dol12]{dolgachev}
Igor~V. Dolgachev.
\newblock {\em Classical Algebraic Geometry: A Modern View}.
\newblock Cambridge University Press, 2012.

\bibitem[Ful04]{fultonblowup}
William Fulton.
\newblock Adjoints and {M}ax {N}oether's {F}undamentalsatz.
\newblock In {\em Algebra, arithmetic and geometry with applications ({W}est
  {L}afayette, {IN}, 2000)}, pages 301--313. Springer, Berlin, 2004.

\bibitem[Ful08]{FultonCurves}
William Fulton.
\newblock {\em Algebraic Curves. An Introduction to Algebraic Geometry}.
\newblock self-published, 2008.

\bibitem[GH81]{grossharris}
Benedict~H. Gross and Joe Harris.
\newblock Real algebraic curves.
\newblock {\em Ann. Sci. \'Ecole Norm. Sup. (4)}, 14(2):157--182, 1981.

\bibitem[GMSY]{RealRootsSource}
Jordy~Lopez Garcia, Kelly Maluccio, Frank Sottile, and Thomas Yahl.
\newblock {RealRoots: A \emph{Macaulay2} package. Version~0.1}.
\newblock A \emph{Macaulay2} package available at
  \url{https://github.com/Macaulay2/M2/tree/master/M2/Macaulay2/packages}.

\bibitem[GMSY24]{M2RealRoots}
Jordy~Lopez Garcia, Kelly Maluccio, Frank Sottile, and Thomas Yahl.
\newblock Real solutions to systems of polynomial equations in {M}acaulay2.
\newblock {\em Journal of Software for Algebra and Geometry}, 14(1):87--95,
  2024.

\bibitem[Gro67]{EGAIV4}
A.~Grothendieck.
\newblock \'{E}l\'ements de g\'eom\'etrie alg\'ebrique. {IV}. \'{E}tude locale
  des sch\'emas et des morphismes de sch\'emas {IV}.
\newblock {\em Inst. Hautes \'Etudes Sci. Publ. Math.}, (32):361, 1967.

\bibitem[GS]{M2}
Daniel~R. Grayson and Michael~E. Stillman.
\newblock Macaulay2, a software system for research in algebraic geometry.
\newblock Available at \url{http://www2.macaulay2.com}.

\bibitem[Har77]{Hart77}
Robin Hartshorne.
\newblock {\em Algebraic geometry}.
\newblock Springer-Verlag, New York-Heidelberg, 1977.
\newblock Graduate Texts in Mathematics, No. 52.

\bibitem[HH22]{harhov}
Davit Harutyunyan and Narek Hovsepyan.
\newblock On the extreme rays of the cone of {$3\times3$} quasiconvex quadratic
  forms: extremal determinants versus extremal and polyconvex forms.
\newblock {\em Arch. Ration. Mech. Anal.}, 244(1):1--25, 2022.

\bibitem[PSV11]{plaumann}
Daniel Plaumann, Bernd Sturmfels, and Cynthia Vincant.
\newblock Quartic curves and their bitangents.
\newblock {\em Journal of Symbolic Computation}, 46:712--733, 2011.

\bibitem[PV13]{hvelemt}
Daniel Plaumann and Cynthia Vinzant.
\newblock Determinantal representations of hyperbolic plane curves: an
  elementary approach.
\newblock {\em J. Symbolic Comput.}, 57:48--60, 2013.

\bibitem[Qua15]{quarez}
Ronan Quarez.
\newblock On the real zeros of positive semidefinite biquadratic forms.
\newblock {\em Comm. Algebra}, 43(3):1317--1353, 2015.

\bibitem[Sch24]{Scheiderer2024RAG}
Claus Scheiderer.
\newblock {\em A course in real algebraic geometry : positivity and sums of
  squares}.
\newblock Number 303 in Graduate Texts in Mathematics (GTM). Springer, Cham,
  2024.

\bibitem[{Sta}24]{stacks-project}
The {Stacks project authors}.
\newblock The stacks project.
\newblock \url{https://stacks.math.columbia.edu}, 2024.

\bibitem[Vin89]{Vinnikov1989LDR}
Victor Vinnikov.
\newblock Complete description of determinantal representations of smooth
  irreducible curves.
\newblock {\em Linear Algebra and its Applications}, 125:103--140, 1989.

\bibitem[Vin93]{vin93}
Victor Vinnikov.
\newblock Selfadjoint determinantal representations of real plane curves.
\newblock {\em Math. Ann.}, 296(3):453--479, 1993.

\end{thebibliography}
